\documentclass{siamltex}
\usepackage{fullpage,enumerate,setspace,changepage,multirow,rotating,color}
\usepackage{graphicx,amsmath,amssymb,url, algorithm, algorithmic}
\def\grad{\nabla}

\def\cK{\mathcal{K}}

\def\cO{\mathcal{O}}

\def\smskip{\smallskip}

\def\texitem#1{\par\smskip\noindent\hangindent 25pt
               \hbox to 25pt {\hss #1 ~}\ignorespaces}

\def\norm#1{\|#1\|}

\newcommand{\BEAS}{\begin{eqnarray*}}
\newcommand{\EEAS}{\end{eqnarray*}}
\newcommand{\BEA}{\begin{eqnarray}}
\newcommand{\EEA}{\end{eqnarray}}
\newcommand{\BEQ}{\begin{eqnarray}}
\newcommand{\EEQ}{\end{eqnarray}}
\newcommand{\BIT}{\begin{itemize}}
\newcommand{\EIT}{\end{itemize}}
\newcommand{\BNUM}{\begin{enumerate}}
\newcommand{\ENUM}{\end{enumerate}}

\newcommand{\BA}{\begin{array}}
\newcommand{\EA}{\end{array}}

\newcommand{\reals}{\mathbb{R}}
\newcommand{\integers}{\mathbb{Z}}

\newcommand{\Rank}{\mathop{\bf rank}}

\newcommand{\argmin}{\mathop{\rm argmin}}
\newcommand{\argmax}{\mathop{\rm argmax}}

\newif\ifpagenumbering
\pagenumberingtrue

\pagenumberingfalse

\newsavebox{\theorembox}
\newsavebox{\lemmabox}
\newsavebox{\remarkbox}
\savebox{\theorembox}{\noindent\bf Theorem}
\savebox{\lemmabox}{\noindent\bf Lemma}
\savebox{\remarkbox}{\noindent\bf Remark}
\newtheorem{remark}{\usebox{\remarkbox}}[section]

\newcommand{\proc}[1]{\textnormal{\scshape#1}}

\def\fprod#1{\langle#1 \rangle}

\title{Fast First-Order Methods for Stable Principal Component Pursuit\footnotemark[5]}
\author{
    N. S. Aybat \footnotemark[2]
    \and
    D. Goldfarb \footnotemark[3]\ 
    \and
    G. Iyengar \footnotemark[4]\ 
}
\begin{document}
\maketitle
\renewcommand{\thefootnote}{\fnsymbol{footnote}}
\footnotetext[2]{IEOR Department, Columbia University. Email: {\tt nsa2106@columbia.edu}.} 
\footnotetext[3]{IEOR Department, Columbia University. Email: {\tt goldfarb@columbia.edu}.}
\footnotetext[4]{IEOR Department, Columbia University. Email: {\tt gi10@columbia.edu}.}
\footnotetext[5]{Research partially supported by ONR grant N000140310514, NSF Grant DMS 10-16571 and DOE Grant DE-FG02-08-25856.}
\renewcommand{\thefootnote}{\arabic{footnote}}
\begin{abstract}
  The stable principal component pursuit~(SPCP) problem is a non-smooth convex optimization problem, the solution of which has been shown both in theory and in practice to enable one to recover the low rank and sparse components of a matrix whose elements have been corrupted by Gaussian noise. In this paper, we first show how several existing fast first-order methods can be applied to this problem very efficiently.  Specifically, we show that the subproblems that arise when applying optimal gradient methods of Nesterov, alternating linearization methods and alternating direction augmented Lagrangian methods to the SPCP problem either have closed-form solutions or have solutions that can be obtained with very modest effort. Later, we develop a new first order algorithm, \proc{NSA}, based on partial variable splitting. All but one of the methods analyzed require at least one of the non-smooth terms in the objective function to be smoothed and obtain an $\epsilon$-optimal solution to the SPCP problem in $O(1/\epsilon)$ iterations. \proc{NSA}, which works directly with the fully non-smooth objective function, is proved to be convergent under mild conditions on the sequence of parameters it uses.  Our preliminary computational tests show that the latter method, \proc{NSA}, although its complexity is not known, is the fastest among the four algorithms described and substantially outperforms \proc{ASALM}, the only existing method for the SPCP problem. To best of our knowledge, an algorithm for the SPCP problem that has $\cO(1/\epsilon)$ iteration complexity and has a per iteration complexity equal to that of a singular value decomposition is given for the first time.
\end{abstract}
\section{Introduction}
In \cite{Can09_1J, Wright09_1J}, it was shown that when the data matrix $D\in\reals^{m\times n}$ is of the form $D = X^0 + S^0$, where $X^0$ is a low-rank matrix, i.e. $\Rank(X^0)\ll\min\{m,n\}$, and $S^0$ is a sparse matrix, i.e. $\norm{S^0}_0\ll mn$ ($\norm{.}_0$ counts the number of nonzero elements of its argument), one can recover the low-rank and sparse components of $D$ by solving the \emph{principal component pursuit} problem
\begin{align}
\min_{X\in\reals^{m\times n}} \norm{X}_* + \xi~\norm{D-X}_1, \label{eq:component_pursuit}
\end{align}
where $\xi=\frac{1}{\sqrt{\max\{m,n\}}}$. 

For $X\in\reals^{m\times n}$, $\norm{X}_*$ denotes the nuclear norm of $X$, which is equal to the sum of its singular values, $\norm{X}_1:=\sum_{i=1}^m\sum_{j=1}^n|X_{ij}|$, $\norm{X}_\infty:=\max\{|X_{ij}|:~1\leq i\leq m,~ 1\leq j\leq n\}$ and $\norm{X}_2:=\sigma_{\rm max}(X)$, where $\sigma_{\rm max}(X)$ is the maximum singular value of $X$.

To be more precise, let $X^0\in\reals^{m\times n}$ with $\Rank(X^0)=r$ and let $X^0=U\Sigma V^T=\sum_{i=1}^r\sigma_iu_iv_i^T$ denote the singular value decomposition~(SVD) of $X^0$. Suppose that for some $\mu>0$, $U$ and $V$ satisfy
\begin{align}
\label{eq:assumption}
\max_i\norm{U^Te_i}_2^2\leq \frac{\mu r}{m}, \quad \max_i\norm{V^Te_i}_2^2\leq \frac{\mu r}{n}, \quad \norm{UV^T}_\infty\leq \sqrt{\frac{\mu r}{mn}},
\end{align}
where $e_i$ denotes the $i$-th unit vector.
\begin{theorem}~\cite{Can09_1J}
Suppose $D=X^0+S^0$, where $X^0\in\reals^{m\times n}$ with $m<n$ satisfies \eqref{eq:assumption} for some $\mu>0$, and the support set of $S^0$ is uniformly distributed. Then there are constants $c$, $\rho_r$, $\rho_s$ such that with probability of at least $1-c n^{-10}$, the principal component pursuit problem~\eqref{eq:component_pursuit} exactly recovers $X^0$ and $S^0$ provided that
\begin{align}
\label{eq:assumption2}
\Rank(X^0)\leq\rho_r m \mu^{-1} (\log(n))^{-2} \quad \mbox{and} \quad \norm{S^0}_0\leq \rho_s mn.
\end{align}
\end{theorem}
In \cite{Can10_1J}, it is shown that the recovery is still possible even when the data matrix, $D$, is corrupted with a dense error matrix, $\zeta^0$ such that $\norm{\zeta^0}_F\leq\delta$, by solving the \emph{stable principal component pursuit}~(SPCP) problem
\begin{align}
(P): \min_{X,S\in\reals^{m\times n}}\{\norm{X}_*+\xi~\norm{S}_1:\ \norm{X+S-D}_F\leq\delta\}. \label{eq:stable_component_pursuit}
\end{align}
Specifically, the following theorem is proved in \cite{Can10_1J}.
\begin{theorem}~\cite{Can10_1J}
\label{thm:candes2}
Suppose $D = X^0 + S^0 + \zeta^0$, where $X^0\in\reals^{m\times n}$ with $m<n$ satisfies \eqref{eq:assumption} for some $\mu>0$, and the support set of $S^0$ is uniformly distributed. If $X^0$ and $S^0$ satisfy \eqref{eq:assumption2}, then for any $\zeta^0$ such that $\norm{\zeta^0}_F\leq\delta$ the solution, $(X^*,S^*)$, to the stable principal component pursuit problem~\eqref{eq:stable_component_pursuit} satisfies $\norm{X^*-X^0}_F^2+\norm{S^*-S^0}_F^2\leq Cmn\delta^2$ for some constant $C$ with high probability.
\end{theorem}

Principal component pursuit and stable principal component pursuit both have applications in video surveillance and face recognition. For existing algorithmic approaches to solving principal component pursuit see ~\cite{Can09_1J,Gold10_1J,Ma09_1J,Ma09_1R,Can10_1J} and references therein. In this paper, we develop four different fast first-order algorithms to solve the 
SPCP problem $(P)$. The first two algorithms are direct applications of Nesterov's optimal algorithm~\cite{Nesterov05} and the proximal gradient method of Tseng~\cite{Tseng08}, which is inspired by both FISTA and Nesterov's infinite memory algorithms that are introduced in \cite{Beck09_1J} and \cite{Nesterov05}, respectively. In this paper it is shown that both algorithms can compute an $\epsilon$-optimal, feasible solution to $(P)$ in $\cO(1/\epsilon)$ iterations. The third and fourth algorithms apply an alternating direction augmented Lagrangian approach to an equivalent problem obtained by partial variable splitting. The third algorithm can compute an $\epsilon$-optimal, feasible solution to the problem in $\cO(1/\epsilon^2)$ iterations, which can be easily improved to $\cO(1/\epsilon)$ complexity. Given $\epsilon>0$, all first three algorithms use suitably smooth versions of at least one of the norms in the objective function. The fourth algorithm~(NSA) works directly with the original non-smooth objective function and can be shown to converge to an optimal solution of $(P)$, provided that a mild condition on the increasing sequence of penalty multipliers holds. To best of our knowledge, an algorithm for the SPCP problem that has $\cO(1/\epsilon)$ iteration complexity and has a per iteration complexity equal to that of a singular value decomposition is given for the first time.

The only algorithm that we know of that has been designed to solve the 
SPCP problem $(P)$ is the algorithm \proc{ASALM}~\cite{Tao09_1J}. The results of our numerical experiments comparing \proc{NSA} algorithm with \proc{ASALM} has shown that \proc{NSA} is faster and also more robust to changes in problem parameters.

\section{Proximal Gradient Algorithm with Smooth Objective Function}
In this section we show that Nesterov's optimal algorithm~\cite{Nesterov04,Nesterov05} for simple sets is efficient for solving $(P)$.

For fixed parameters $\mu>0$ and $\nu>0$, define the smooth $C^{1,1}$ functions $f_\mu(.)$ and $g_\nu(.)$ as follows
\begin{align}
&f_\mu(X)=\max_{U\in\reals^{m\times n}:\norm{U}_2\leq 1}\fprod{X,U}-\frac{\mu}{2}\norm{U}_F^2, \label{eq:smooth_f}\\
&g_\nu(S)=\max_{W\in\reals^{m\times n}:\norm{W}_\infty\leq 1}\fprod{S,W}-\frac{\nu}{2}\norm{W}_F^2. \label{eq:smooth_g}
\end{align}
Clearly, $f_\mu(.)$ and $g_\nu(.)$ closely approximate the non-smooth functions $f(X):=\norm{X}_*$ and $g(S):=\norm{S}_1$, respectively. Also let $\chi:=\{(X,S)\in\reals^{m\times n}\times\reals^{m\times n}:~\norm{X+S-D}_F\leq\delta\}$ and $L=\frac{1}{\mu}+\frac{1}{\nu}$, where $\frac{1}{\mu}$ and $\frac{1}{\nu}$ are the Lipschitz constants for the gradients of $f_\mu(.)$ and $g_\nu(.)$, respectively. Then Nesterov's optimal algorithm~\cite{Nesterov04,Nesterov05} for simple sets applied to the problem:
\begin{align}
\min_{X,S\in\reals^{m\times n}}\{f_\mu(X)+\xi~g_\nu(S):\ (X,S)\in\chi\}, \label{eq:problem_ss}
\end{align}
is given by \textbf{Algorithm~\ref{alg:sNesterov}}.
\begin{algorithm}[h!]
    \caption{SMOOTH PROXIMAL GRADIENT($X_0,S_0$)}\label{alg:sNesterov}
    {\small
    \begin{algorithmic}[1]
    \STATE \textbf{input:} $X_0\in\reals^{m\times n}$, $S_0\in\reals^{m\times n}$
    \STATE $k\gets 0$
    \WHILE{$k\leq k^*$}
    \STATE Compute $\grad f_\mu(X_k)$ and $\grad g_\nu(S_k)$
    \STATE $(Y^x_k, Y^s_k)\gets\argmin_{X,S}\left\{\fprod{\grad f_\mu(X_k), X}+\fprod{\grad g_\nu(S_k), S}+\frac{L}{2}\left(\norm{X-X_k}_F^2+\norm{S-S_k}_F^2\right):\ (X,S)\in \chi\right\}$
    \STATE $\Gamma_k(X,S):=\sum_{i=0}^k \frac{i+1}{2}\left\{\fprod{\grad f_\mu(X_i), X}+\fprod{\grad g_\nu(S_i), S}\right\}$
    \STATE $(Z^x_k, Z^s_k)\gets\argmin_{X,S}\left\{\Gamma_k(X,S)+ \frac{L}{2}\left(\norm{X-X_0}_F^2+\norm{S-S_0}_F^2\right):(X,S)\in \chi\right\}$
    \STATE $(X_{k+1},S_{k+1})\gets \left(\frac{k+1}{k+3}\right)(Y^x_k, Y^s_k)+\left(\frac{2}{k+3}\right)(Z^x_k, Z^s_k)$
    \STATE $k \gets k + 1$
    \ENDWHILE
    \RETURN $(X_{k^*},S_{k^*})$
    \end{algorithmic}
    }
\end{algorithm}

Because of the simple form of the set $\chi$, it is easy to ensure that all iterates $(Y^x_k, Y^s_k)$, $(Z^x_k, Z^s_k)$ and $(X_{k+1},S_{k+1})$ lie in $\chi$.
Hence, \textbf{Algorithm~\ref{alg:sNesterov}} enjoys the full convergence rate of $\cO(L/k^2)$ of the Nesterov's method. Thus, setting $\mu=\Omega(\epsilon)$ and $\nu=\Omega(\epsilon)$, \textbf{Algorithm~\ref{alg:sNesterov}} computes an $\epsilon$-optimal and feasible solution to problem $(P)$ in $k^*=\cO(1/\epsilon)$ iterations.
The iterates $(Y^x_k, Y^s_k)$ and $(Z^x_k, Z^s_k)$ that need to be computed at each iteration of \textbf{Algorithm~\ref{alg:sNesterov}} are solutions to an optimization problem of the form:
\begin{align}
\label{eq:subproblem_sa}
(P_s):\  \min_{X,S\in\reals^{m\times n}}\left\{\frac{L}{2}\left(\norm{X-\tilde{X}}_F^2+\norm{S-\tilde{S}}_F^2\right)+\fprod{Q_x, X}+\fprod{Q_s, S}:\ (X,S)\in\chi\right\}.
\end{align}
The following lemma shows that the solution to problems of the form $(P_s)$ can be computed efficiently.
\newpage
\begin{lemma}
\label{lem:subproblem_smoothsmooth}
The optimal solution $(X^*,S^*)$ to problem $(P_s)$
can be written in closed form as follows.

When $\delta>0$,
\begin{align}
&X^*=\left(\frac{\theta^*}{L+2\theta^*}\right)\left(D-q_s(\tilde{S})\right)+\left(\frac{L+\theta^*}{L+2\theta^*}\right)q_x(\tilde{X}), \label{lemeq:X_smooth}\\
&S^*=\left(\frac{\theta^*}{L+2\theta^*}\right)\left(D-q_x(\tilde{X})\right)+\left(\frac{L+\theta^*}{L+2\theta^*}\right)q_s(\tilde{S}), \label{lemeq:S_smooth}
\end{align}
where $q_x(X):=X-\frac{1}{L}~Q_x$, $q_s(S):=S-\frac{1}{L}~Q_s$ and
\begin{align}
\theta^*=\max\left\{0,~\frac{L}{2}\left(\frac{\norm{q_x(\tilde{X})+q_s(\tilde{S})-D}_F}{\delta}-1\right)\right\}.
\end{align}

When $\delta=0$,
\begin{equation}
\label{lemeq:XS_smooth_delta0}
\begin{array}{ll}
X^*=\frac{1}{2}\left(D-q_s(\tilde{S})\right)+\frac{1}{2}~q_x(\tilde{X}) \hbox{ and} &S^*=\frac{1}{2}\left(D-q_x(\tilde{X})\right)+\frac{1}{2}~q_s(\tilde{S}).
\end{array}
\end{equation}
\end{lemma}
\begin{proof}
Suppose that $\delta>0$. Writing the constraint in problem $(P_s)$, $(X,S)\in\chi$, as
\begin{align}
\frac{1}{2}\norm{X+S-D}^2_F\leq \frac{\delta^2}{2}, \label{eq:quadratic_constraint}
\end{align}
the Lagrangian function for \eqref{eq:subproblem_sa} is given as
\begin{align*}
\mathcal{L}(X,S;\theta)=\frac{L}{2}\left(\norm{X-\tilde{X}}_F^2+\norm{S-\tilde{S}}_F^2\right)+\fprod{Q_x, X-\tilde{X}}+\fprod{Q_s, S-\tilde{S}}+\frac{\theta}{2}\left(\norm{X+S-D}_F^2-\delta^2\right).
\end{align*}
Therefore, the optimal solution $(X^*,S^*)$ and optimal Lagrangian multiplier $\theta^*\in\reals$ must satisfy the Karush-Kuhn-Tucker~(KKT) conditions:
\begin{enumerate}[i.]
    \item $\norm{X^*+S^*-D}_F\leq\delta$, \label{condition1_sa}
    \item $\theta^*\geq 0$, \label{condition2_sa}
    \item $\theta^*~(\norm{X^*+S^*-D}_F-\delta)=0$, \label{condition3_sa}
    \item $L(X^*-\tilde{X})+\theta^*(X^*+S^*-D)+Q_x=0$, \label{condition4_sa}
    \item $L(S^*-\tilde{S})+\theta^*(X^*+S^*-D)+Q_s=0$. \label{condition5_sa}
\end{enumerate}
Conditions \ref{condition4_sa} and \ref{condition5_sa} imply that $(X^*,S^*)$ satisfy $\eqref{lemeq:X_smooth}$ and \eqref{lemeq:S_smooth}, from which it follows that
\begin{align}
X^*+S^*-D=\left(\frac{L}{L+2\theta^*}\right)\left(q_x(\tilde{X})+q_s(\tilde{S})-D\right). \label{eq:smooth_infeasibility}
\end{align}
\subsection*{Case 1: $\norm{q_x(\tilde{X})+q_s(\tilde{S})-D}_F\leq\delta$}
Setting $X^*=q_x(\tilde{X})$, $S^*=q_s(\tilde{S})$ and $\theta^*=0$, clearly satisfies $\eqref{lemeq:X_smooth}$, \eqref{lemeq:S_smooth} and conditions \ref{condition1_sa} (from \eqref{eq:smooth_infeasibility}), \ref{condition2_sa} and \ref{condition3_sa}. Thus, this choice of variables satisfies all the five KKT conditions.
\subsection*{Case 2: $\norm{q_x(\tilde{X})+q_s(\tilde{S})-D}_F>\delta$}
Set $\theta^*=\frac{L}{2}\left(\frac{\norm{q_x(\tilde{X})+q_s(\tilde{S})-D}_F}{\delta}-1\right)$. Since $\norm{q_x(\tilde{X})+q_s(\tilde{S})-D}_F>\delta$, $\theta^*>0$; hence, \ref{condition2_sa} is satisfied. Moreover, for this value of $\theta^*$, it follows from \eqref{eq:smooth_infeasibility} that $\norm{X^*+S^*-D}_F=\delta$. Thus, KKT conditions \ref{condition1_sa} and \ref{condition3_sa} are satisfied.

Therefore, setting $X^*$ and $S^*$ according to \eqref{lemeq:X_smooth} and \eqref{lemeq:S_smooth}, respectively; and setting
\begin{align*}
\theta^*= \max\left\{0,\ \frac{L}{2}\left(\frac{\norm{q_x(\tilde{X})+q_s(\tilde{S})-D}_F}{\delta}-1\right)\right\},
\end{align*}
satisfies all the five KKT conditions.

Now, suppose that $\delta=0$. Since $S^*=D-X^*$, problem~$(P_s)$
can be written as
\begin{equation*}
\begin{array}{ll}
\min_{X\in\reals^{m\times n}} & \norm{X-\tilde{X}+\frac{Q_x}{L}}_F^2+\norm{D-X-\tilde{S}+\frac{Q_s}{L}}_F^2,
\end{array}
\end{equation*}
which is also equivalent to the problem: $\min_{X\in\reals^{m\times n}} \norm{X-q_x(\tilde{X})}_F^2+\norm{X-(D-q_s(\tilde{S}))}_F^2$. Then \eqref{lemeq:XS_smooth_delta0} trivially follows from first-order optimality conditions for this problem and the fact that $S^*=D-X^*$.
\end{proof}

\section{Proximal Gradient Algorithm with Partially Smooth Objective Function}
In this section we show how the proximal gradient algorithm, Algorithm~3 in \cite{Tseng08}, can be applied to the problem
\begin{align}
\label{eq:problem_sns}
\min_{X,S\in\reals^{m\times n}}\{f_\mu(X)+\xi~\norm{S}_1:\ (X,S)\in\chi\},
\end{align}
where $f_\mu(.)$ is the smooth function defined in \eqref{eq:smooth_f} such that $\grad f_\mu(.)$ is Lipschitz continuous with constant $L_\mu=\frac{1}{\mu}$.
This algorithm is given in \textbf{Algorithm~\ref{alg:nsNesterov}}.
\begin{algorithm}[h!]
    \caption{PARTIALLY SMOOTH PROXIMAL GRADIENT($X_0,S_0$)}\label{alg:nsNesterov}
    {\small
    \begin{algorithmic}[1]
    \STATE \textbf{input:} $X_0\in\reals^{m\times n}$, $S_0\in\reals^{m\times n}$
    \STATE $(Z^x_0, Z^s_0)\gets(X_0,S_0)$, $k\gets 0$
    \WHILE{$k\leq k^*$}
    \STATE $(Y^x_k, Y^s_k)\gets \left(\frac{k}{k+2}\right)(X_{k},S_{k})+\left(\frac{2}{k+2}\right)(Z^x_k, Z^s_k)$
    \STATE Compute $\grad f_\mu(Y^x_k)$
    \STATE $(Z^x_{k+1}, Z^s_{k+1})\gets\argmin_{X,S}\left\{\sum_{i=0}^k \frac{i+1}{2}\left\{\xi\norm{S}_1+\fprod{\grad f_\mu(Y^x_i), X}\right\}+\frac{L_\mu}{2}\norm{X-X_0}_F^2:\ (X,S)\in \chi\right\}$
    \STATE $(X_{k+1},S_{k+1})\gets \left(\frac{k}{k+2}\right)(X_{k},S_{k})+\left(\frac{2}{k+2}\right)(Z^x_{k+1}, Z^s_{k+1})$
    \STATE $k \gets k + 1$
    \ENDWHILE
    \RETURN $(X_{k^*},S_{k^*})$
    \end{algorithmic}
    }
\end{algorithm}

Mimicking the proof in \cite{Tseng08}, it is easy to show that \textbf{Algorithm~\ref{alg:nsNesterov}}, which uses the prox function $\frac{1}{2}\norm{X-X_0}_F^2$, converges to the optimal solution of \eqref{eq:problem_sns}. Given $(X_0,S_0)\in\chi$, e.g. $X_0=\mathbf{0}$ and $S_0=D$, the current algorithm keeps all iterates in $\chi$ as in \textbf{Algorithm~\ref{alg:sNesterov}}, and hence it enjoys the full convergence rate of $\cO(L/k^2)$. Thus, setting $\mu=\Omega(\epsilon)$, \textbf{Algorithm~\ref{alg:nsNesterov}} computes an $\epsilon$-optimal, feasible solution of problem $(P)$ in $k^*=\cO(1/\epsilon)$ iterations.

The only thing left to be shown is that the optimization subproblems in \textbf{Algorithm~\ref{alg:nsNesterov}} can be solved efficiently. The subproblem that has to be solved at each iteration to compute $(Z^x_{k+1}, Z^s_{k+1})$ has the form:
\begin{align}
\label{eq:subproblem_nsa}
(P_{ns}):\  \min\left\{\xi\norm{S}_1+\fprod{Q, X-\tilde{X}}+\frac{\rho}{2}\norm{X-\tilde{X}}_F^2:\ (X,S)\in \chi\right\},
\end{align}
for some $\rho>0$. Lemma~\ref{lem:subproblem} shows that these computations can be done efficiently.
\newpage
\begin{lemma}
\label{lem:subproblem}
The optimal solution $(X^*,S^*)$ to problem $(P_{ns})$
can be written in closed form as follows.

When $\delta>0$,
\begin{align}
&S^*=sign\left(D-q(\tilde{X})\right)\odot\max\left\{|D-q(\tilde{X})|-\xi\frac{(\rho+\theta^*)}{\rho\theta^*}~E,\ \mathbf{0}\right\}, \label{lemeq:S}\\
&X^*= \frac{\theta^*}{\rho+\theta^*}~(D-S^*)+\frac{\rho}{\rho+\theta^*}~q(\tilde{X}), \label{lemeq:X}
\end{align}
where $q(\tilde{X}):=\tilde{X}-\frac{1}{\rho}~Q$, $E$ and $\mathbf{0}\in\reals^{m\times n}$ are matrices with all components equal to ones and zeros, respectively, and $\odot$ denotes the componentwise multiplication operator. $\theta^*=0$ if $\norm{D-q(\tilde{X})}_F\leq\delta$; otherwise, $\theta^*$ is the unique positive solution of the nonlinear equation $\phi(\theta)=\delta$, where
\begin{align}
\phi(\theta):= \norm{\min\left\{\frac{\xi}{\theta}~E,\ \frac{\rho}{\rho+\theta}~|D-q(\tilde{X})|\right\}}_F.
\end{align}
Moreover, $\theta^*$ can be efficiently computed in $\cO(mn\log(mn))$ time.

When $\delta=0$,
\begin{equation}
\label{lemeq:XS_nonsmooth_delta0}
\begin{array}{cc}
S^*=sign\left(D-q(\tilde{X})\right)\odot\max\left\{|D-q(\tilde{X})|-\frac{\xi}{\rho}~E,\ \mathbf{0}\right\} \hbox{ and } &X^*=D-S^*.
\end{array}
\end{equation}
\end{lemma}
\begin{proof}
Suppose that $\delta>0$. Let $(X^*,S^*)$ be an optimal solution to problem $(P_{ns})$ 
and $\theta^*$ denote the optimal Lagrangian multiplier for the constraint $(X,S)\in\chi$ written as \eqref{eq:quadratic_constraint}. Then the KKT optimality conditions for this problem are
\begin{enumerate}[i.]
    \item $Q+\rho(X^*-\tilde{X})+\theta^*(X^*+S^*-D)=0$, \label{condition1}
    \item $\xi G + \theta^*(X^*+S^*-D)=0$ and $G\in\partial\norm{S^*}_1$, \label{condition2}
    \item $\norm{X^*+S^*-D}_F\leq\delta$, \label{condition3}
    \item $\theta^*\geq 0$, \label{condition4}
    \item $\theta^*~(\norm{X^*+S^*-D}_F-\delta)=0$. \label{condition5}
\end{enumerate}

From \ref{condition1} and \ref{condition2}, we have
\begin{eqnarray}
\left[
  \begin{array}{cc}
    (\rho+\theta^*)I &  \theta^*I\\
    \theta^*I & \theta^*I \\
  \end{array}
\right]
\left[
  \begin{array}{c}
    X^* \\
    S^* \\
  \end{array}
\right]
=
\left[
  \begin{array}{c}
    \theta^*D+\rho~q(\tilde{X}) \\
    \theta^*D-\xi G\\
  \end{array}
\right], \label{eq:FTOC_1}
\end{eqnarray}
where $q(\tilde{X})=\tilde{X}-\frac{1}{\rho}~Q$. From \eqref{eq:FTOC_1} it follows that
\begin{eqnarray}
\left[
  \begin{array}{cc}
    (\rho+\theta^*)I &  \theta^*I\\
    0 & \left(\frac{\rho\theta^*}{\rho+\theta^*}\right)~I \\
  \end{array}
\right]
\left[
  \begin{array}{c}
    X^* \\
    S^* \\
  \end{array}
\right]
=
\left[
  \begin{array}{c}
    \theta^*D+\rho~q(\tilde{X}) \\
    \frac{\rho\theta^*}{\rho+\theta^*}~(D-q(\tilde{X}))-\xi G\\
  \end{array}
\right]. \label{eq:FTOC_2}
\end{eqnarray}
From the second equation in \eqref{eq:FTOC_2}, we have
\begin{align}
\xi\frac{(\rho+\theta^*)}{\rho\theta^*}~G+S^*+q(\tilde{X})-D=0. \label{eq:shrinkS}
\end{align}
But \eqref{eq:shrinkS} is precisely the first-order optimality conditions for the ``shrinkage" problem $$\min_{S\in\reals^{m\times  n}}\left\{\xi\frac{(\rho+\theta^*)}{\rho\theta^*}\norm{S}_1+\frac{1}{2}\norm{S+q(\tilde{X})-D}_F^2\right\}.$$ Thus, $S^*$ is the optimal solution to the ``shrinkage" problem and is given by \eqref{lemeq:S}.
\eqref{lemeq:X} follows from the first equation in \eqref{eq:FTOC_2}, and it implies
\begin{align}
X^*+S^*-D = \frac{\rho}{\rho+\theta^*}~(S^*+q(\tilde{X})-D).
\end{align}
Therefore,
\begin{align}
\norm{X^*+S^*-D}_F &= \frac{\rho}{\rho+\theta^*}~\norm{S^*+q(\tilde{X})-D}_F, \nonumber\\
&=\frac{\rho}{\rho+\theta^*}~\norm{sign\left(D-q(\tilde{X})\right)\odot\max\left\{|D-q(\tilde{X})|-\xi\frac{(\rho+\theta^*)}{\rho\theta^*}~E,\ \mathbf{0}\right\}-\left(D-q(\tilde{X})\right)}_F, \nonumber\\
&=\frac{\rho}{\rho+\theta^*}~\norm{\max\left\{|D-q(\tilde{X})|-\xi\frac{(\rho+\theta^*)}{\rho\theta^*}~E,\ \mathbf{0}\right\}-|D-q(\tilde{X})|~}_F,\nonumber\\
&=\frac{\rho}{\rho+\theta^*}~\norm{\min\left\{\xi\frac{(\rho+\theta^*)}{\rho\theta^*}~E,\ |D-q(\tilde{X})|\right\}}_F,\nonumber\\
&=\norm{\min\left\{\frac{\xi}{\theta^*}~E,\ \frac{\rho}{\rho+\theta^*}~|D-q(\tilde{X})|\right\}}_F, \label{eq:Fnorm}
\end{align}
where the second equation uses \eqref{lemeq:S}. Now let $\phi:\reals_+\rightarrow\reals_+$ be
\begin{align}
\phi(\theta):= \norm{\min\left\{\frac{\xi}{\theta}~E,\ \frac{\rho}{\rho+\theta}~|D-q(\tilde{X})|\right\}}_F.
\end{align}
\subsection*{Case 1: $\norm{D-q(\tilde{X})}_F\leq\delta$}
$\theta^*=0$, $S^*=\mathbf{0}$ and $X^*=q(\tilde{X})$ trivially satisfy all the KKT conditions.
\subsection*{Case 2: $\norm{D-q(\tilde{X})}_F>\delta$}
It is easy to show that $\phi(.)$ is a strictly decreasing function of $\theta$. Since $\phi(0)=\norm{D-q(\tilde{X})}_F>\delta$ and $\lim_{\theta\rightarrow\infty}\phi(\theta)=0$, there exists a unique $\theta^*>0$ such that $\phi(\theta^*)=\delta$. Given $\theta^*$, $S^*$ and $X^*$ can then be computed from equations \eqref{lemeq:S} and \eqref{lemeq:X}, respectively. Moreover, since $\theta^*>0$ and $\phi(\theta^*)=\delta$, \eqref{eq:Fnorm} implies that $X^*$, $S^*$ and $\theta^*$ satisfy the KKT conditions.

We now show that $\theta^*$ can be computed in $\cO(mn\log(mn))$ time. Let $A:=|D-q(\tilde{X})|$ and $0\leq a_{(1)}\leq a_{(2)}\leq ... \leq a_{(mn)}$ be the $mn$ elements of the matrix $A$ sorted in increasing order, which can be done in $\cO(mn\log(mn))$ time. Defining $a_{(0)}:=0$ and $a_{(mn+1)}:=\infty$, we then have for all $j\in\{0,1,...,mn\}$ that
\begin{align}
\frac{\rho}{\rho+\theta}~a_{(j)} \leq \frac{\xi}{\theta} \leq \frac{\rho}{\rho+\theta}~a_{(j+1)} \Leftrightarrow \frac{1}{\xi}~a_{(j)}-\frac{1}{\rho} \leq \frac{1}{\theta} \leq \frac{1}{\xi}~a_{(j+1)}-\frac{1}{\rho}.
\end{align}
For all $\bar{k}< j\leq mn$ define $\theta_j$ such that $\frac{1}{\theta_j}=\frac{1}{\xi}~a_{(j)}-\frac{1}{\rho}$ and let $\bar{k}:=\max\left\{j: \frac{1}{\theta_j}\leq 0,\ j\in\{0,1,...,mn\}\right\}$. Then for all $\bar{k}< j\leq mn$
\begin{align}
\phi(\theta_j)=\sqrt{\left(\frac{\rho}{\rho+\theta_j}\right)^2~\sum_{i=0}^j a^2_{(i)}+(mn-j)~\left(\frac{\xi}{\theta_j}\right)^2}.
\end{align}
Also define $\theta_{\bar{k}}:=\infty$ and $\theta_{mn+1}:=0$ so that $\phi(\theta_{\bar{k}}):=0$ and $\phi(\theta_{mn+1})=\phi(0)=\norm{A}_F>\delta$. Note that $\{\theta_j\}_{\{\bar{k}< j\leq mn\}}$ contains all the points at which $\phi(\theta)$ may not be differentiable for $\theta\geq 0$.
Define $j^*:=\max\{j:\ \phi(\theta_j)\leq\delta,\ \bar{k}\leq j\leq mn\}$. Then $\theta^*$ is the unique solution of the system
\begin{align}
\label{eq:root}
\sqrt{\left(\frac{\rho}{\rho+\theta}\right)^2~\sum_{i=0}^{j^*} a^2_{(i)}+(mn-j^*)~\left(\frac{\xi}{\theta}\right)^2}=\delta \mbox{ and } \theta>0,
\end{align}
since $\phi(\theta)$ is continuous and strictly decreasing in $\theta$ for $\theta\geq 0$. Solving the equation in \eqref{eq:root} requires finding the roots of a fourth-order polynomial (a.k.a. quartic function); therefore, one can compute $\theta^*>0$ using the algebraic solutions of quartic equations (as shown by Lodovico Ferrari in 1540), which requires $\cO(1)$ operations.

Note that if $\bar{k}=mn$, then $\theta^*$ is the solution of the equation
\begin{align}
\sqrt{\left(\frac{\rho}{\rho+\theta^*}\right)^2~\sum_{i=1}^{mn} a^2_{(i)}}=\delta,
\end{align}
i.e. $\theta^*= \rho\left(\frac{\norm{A}_F}{\delta}-1\right)=\rho\left(\frac{\norm{D-\tilde{X}}_F}{\delta}-1\right)$.
Hence, we have proved that problem~$(P_{ns})$ can be solved efficiently.

Now, suppose that $\delta=0$. Since $S^*=D-X^*$, problem~$(P_{ns})$ 
can be written as
\begin{equation}
\label{eq:subproblem_delta0}
\begin{array}{ll}
\min_{S\in\reals^{m\times n}} & \frac{\xi}{\rho}\norm{S}_1+\frac{1}{2}\norm{S-(D-q(\tilde{X}))}_F^2.
\end{array}
\end{equation}
Then \eqref{lemeq:XS_nonsmooth_delta0} trivially follows from first-order optimality conditions for the above problem and the fact that $X^*=D-S^*$.
\end{proof}

The following lemma will be used later in Section~\ref{sec:nsa}. However, we give its proof here, since it uses some equations from the proof of Lemma~\ref{lem:subproblem}. Let $\mathbf{1}_\chi(.,.)$ denote the indicator function of the closed convex set $\chi\subset\reals^{m\times n}\times\reals^{m\times n}$, i.e. if $(Z,S)\in\chi$, then $\mathbf{1}_\chi(Z,S)=0$; otherwise, $\mathbf{1}_\chi(Z,S)=\infty$.
\begin{lemma}
\label{lem:chi_subgradient}
Suppose that $\delta>0$. Let $(X^*,S^*)$ be an optimal solution to problem~$(P_{ns})$
and $\theta^*$ be an optimal Lagrangian multiplier such that $(X^*,S^*)$ and $\theta^*$ together satisfy the KKT conditions, \ref{condition1}-\ref{condition5} in the proof of Lemma~\ref{lem:subproblem}. Then $(W^*,W^*)\in\partial \mathbf{1}_\chi(X^*,S^*)$, where 
$W^*:=-Q+\rho(\tilde{X}-X^*)=\theta^*(X^*+S^*-D)$.
\end{lemma}
\begin{proof}
Let $W^*:=-Q+\rho(\tilde{X}-X^*)$, then from \ref{condition1} and \ref{condition5} of the KKT optimality conditions
in the proof of Lemma~\ref{lem:subproblem}, we have $W^*=\theta^*(X^*+S^*-D)$ and
\begin{align}
\norm{W^*}_F=\theta^*\norm{X^*+S^*-D}=\theta^*(\norm{X^*+S^*-D}-\delta)+\theta^*\delta=\theta^*\delta.
\end{align}
Moreover, for all $(X,S)\in\chi$, it follows from the definition of $\chi$ that $\fprod{W^*, \theta^*(X+S-D)}\leq\theta^*\norm{W^*}_F\norm{X+S-D}_F\leq\theta^*\delta\norm{W^*}_F$. Thus, for all $(X,S)\in\chi$, we have $\fprod{W^*,W^*}=\norm{W^*}_F^2=\theta^*\delta\norm{W^*}_F\geq\fprod{W^*, \theta^*(X+S-D)}$.
Hence,
\begin{align}
0\geq\fprod{W^*, \theta^*(X+S-D)-W^*}=\fprod{W^*, \theta^*(X-X^*+S-S^*)}\hspace{5mm}\forall~(X,S)\in\chi. \label{eq:subgradient_key}
\end{align}
It follows from the proof of Lemma~\ref{lem:subproblem} that if $\norm{D-q(\tilde{X})}_F>\delta$, then $\theta^*>0$, where $q(\tilde{X})=\tilde{X}-\frac{1}{\rho}Q$. Therefore, \eqref{eq:subgradient_key} implies that
\begin{align}
0\geq\fprod{W^*, X-X^*+S-S^*}\hspace{5mm}\forall~(X,S)\in\chi. \label{eq:subgradient_ineq}
\end{align}
On the other hand, if $\norm{D-q(\tilde{X})}_F\leq\delta$, then $\theta^*=0$. Hence $W^*=\theta^*(X^*+S^*-D)=0$, and \eqref{eq:subgradient_ineq} follows trivially. Therefore, \eqref{eq:subgradient_ineq} always holds and this shows that $(W^*,W^*)\in\partial \mathbf{1}_\chi(X^*,S^*)$.
\end{proof}
\section{Alternating Linearization and Augmented Lagrangian Algorithms}
In this and the next section we present algorithms for solving problems 
\eqref{eq:problem_sns} and \eqref{eq:stable_component_pursuit} that are based on partial variable splitting combined with alternating minimization of a suitably linearized augmented Lagrangian function.
We can  write problems 
\eqref{eq:stable_component_pursuit} and \eqref{eq:problem_sns} generically as
\begin{align}
\label{eq:problem_generic}
\min_{X,S\in\reals^{m\times n}}\{\phi(X)+\xi~g(S):\ (X,S)\in\chi\}.
\end{align}
For problem~\eqref{eq:stable_component_pursuit}, $\phi(X) = f(X) = \norm{X}_{\ast}$, while for problem~\eqref{eq:problem_sns}, $\phi(X) = f_{\mu}(X)$ given in \eqref{eq:smooth_f}.

In this section, we first assume that assume that $\phi:\reals^{m\times n}\rightarrow \reals$ and $g:\reals^{m\times n}\times\reals^{m\times n}\rightarrow\reals$ are any closed convex functions such that $\grad \phi$ is Lipschitz continuous, and  $\chi$ is a general closed convex set.
Here we use partial variable splitting, i.e. we only split the $X$ variables in \eqref{eq:problem_generic}, to arrive at the following equivalent problem
\begin{align}
\label{eq:problem_generic_split}
\min_{X,S,Z\in\reals^{m\times n}}\{\phi(X)+\xi~g(S):\ X=Z,\ (Z,S)\in\chi\}.
\end{align}
Let $\psi(Z,S):=\xi~g(S)+\mathbf{1}_\chi(Z,S)$ and define the augmented Lagrangian function
\begin{align}
\label{eq:augmemted_lagrangian}
\mathcal{L}_\rho(X,Z,S;Y)=\phi(X)+\psi(Z,S)+\fprod{Y, X-Z}+\frac{\rho}{2}\norm{X-Z}_F^2.
\end{align}
Then minimizing \eqref{eq:augmemted_lagrangian} by alternating between  $X$ and then $(Z,S)$
leads to several possible methods that can compute a solution to \eqref{eq:problem_generic_split}. These include the alternating linearization method~(ALM) 
with skipping step 
that has an
$\cO(\frac{\rho}{k})$ convergence rate, and the fast version 
of this method with an $\cO(\frac{\rho}{k^2})$ rate (see~\cite{Gold10_1J} for full splitting versions of these methods). In this paper, we only provide a proof of the complexity result for the alternating linearization method with skipping steps~(\proc{ALM-S}) in Theorem~\ref{thm:alm} 
below. One can easily extend the proof of Theorem~\ref{thm:alm} to an ALM method based on \eqref{eq:augmemted_lagrangian} with the function $g(S)$ replaced by a suitably smoothed version (see~\cite{Gold10_1J} for the details of ALM algorithm).

\begin{algorithm}[h!]
    \caption{ALM-S($Y_0$)}\label{alg:alms}
    {\small
    \begin{algorithmic}[1]
    \STATE \textbf{input:} $X_0\in\reals^{m\times n}$, $S_0\in\reals^{m\times n}$, $Y_0\in\reals^{m\times n}$
    \STATE $Z_0\gets X_0$, $k \gets 0$
    \WHILE{$k\geq 0$}
    \STATE $X_{k+1}\gets\argmin_X \mathcal{L}_\rho(X,Z_k,S_k;Y_k)$ \label{algeq:alms_subproblem1}
    \IF    {$\phi(X_{k+1})+\psi(X_{k+1},S_k)>\mathcal{L}_\rho(X_{k+1},Z_k,S_k;Y_k)$}
    \STATE $X_{k+1}\gets Z_k$
    \ENDIF
    \STATE $(Z_{k+1},S_{k+1})\gets\argmin_{Z,S} \psi(Z,S)+\phi(X_{k+1})+\fprod{\grad \phi(X_{k+1}), Z-X_{k+1}}+\frac{\rho}{2}\norm{Z-X_{k+1}}_F^2$ \label{algeq:alms_subproblem2}
    \STATE $Y_{k+1}\gets -\grad \phi(X^{k+1})+\rho (X_{k+1}-Z_{k+1})$
    \STATE $k \gets k + 1$
    \ENDWHILE
    \end{algorithmic}
    }
\end{algorithm}

\begin{theorem}
\label{thm:alm}
Let $\phi:\reals^{m\times n}\rightarrow \reals$ and $\psi:\reals^{m\times n}\times\reals^{m\times n}\rightarrow\reals$ be closed convex functions such that $\grad \phi$ is Lipschitz continuous with Lipschitz constant $L$, and $\chi$ be a closed convex set. Let $\Phi(X,S):=\phi(X)+\psi(X,S)$. For $\rho\geq L$, the sequence $\{Z_k,S_k\}_{k\in\integers_+}$ in Algorithm \proc{ALM-S} satisfies
\begin{align}
\label{eq:alms_thm}
\Phi(Z_k,S_k)-\Phi(X^*,S^*)\leq\rho~\frac{\norm{X_0-X^*}_F^2}{2(k+n_k)},
\end{align}
where $(X^*,S^*)=\argmin_{X,S\in\reals^{m\times n}}\Phi(X,S)$, $n_k:=\sum_{i=0}^{k-1}\mathbf{1}_{\{\Phi(X_{i+1},S_i)>\mathcal{L}_\rho(X_{i+1},Z_i,S_i;Y_i)\}}$ and $\mathbf{1}_{\{.\}}$ is 1 if its argument is true; otherwise, 0.
\end{theorem}
\begin{proof}
See Appendix~\ref{sec:thm_proof} for the proof.
\end{proof}

We obtain \textbf{Algorithm~\ref{alg:sns_alm}} by applying \textbf{Algorithm~\ref{alg:alms}} to solve problem \eqref{eq:problem_sns}, where the smooth function $ \phi(X) = f_\mu(X)$, defined in \eqref{eq:smooth_f}, the non-smooth closed convex function is $\xi~\norm{S}_1+\mathbf{1}_{\chi}(X,S)$ and $\chi=\{(X,S)\in\reals^{m\times n}\times\reals^{m\times n}:\ \norm{X+S-D}_F\leq\delta\}$. Theorem~\ref{thm:alm} shows that \textbf{Algorithm~\ref{alg:sns_alm}} has an  iteration complexity of $\cO(\frac{1}{\epsilon^2})$ to obtain $\epsilon$-optimal and feasible solution of $(P)$.
\begin{algorithm}[h!]
    \caption{PARTIALLY SMOOTH ALM($Y_0$)}\label{alg:sns_alm}
    {\small
    \begin{algorithmic}[1]
    \STATE \textbf{input:} $Y_0\in\reals^{m\times n}$
    \STATE $Z_0\gets 0$, $S_0\gets D$, $k \gets 0$
    \WHILE{$k\geq 0$}
    \STATE $X_{k+1}\gets\argmin_X f_\mu(X)+\fprod{Y_k, X-Z_k}+\frac{\rho}{2}\norm{X-Z_k}_F^2$ \label{eq:sns_alm_subproblem1}
    \STATE $B_k\gets f_\mu(X_{k+1})+\xi~\norm{S_k}_1+\fprod{Y_k, X_{k+1}-Z_k}+\frac{\rho}{2}\norm{X_{k+1}-Z_k}_F^2$
    \IF    {$f_\mu(X_{k+1})+\xi~\norm{S_k}_1+\mathbf{1}_{\chi}(X_{k+1},S_k)> B_k$}
    \STATE $X_{k+1}\gets Z_k$
    \ENDIF
    \STATE $(Z_{k+1},S_{k+1})\gets\argmin_{Z,S}\{\xi~\norm{S}_1+\fprod{\grad f_\mu(X_{k+1}), Z-X_{k+1}}+\frac{\rho}{2}\norm{Z-X_{k+1}}_F^2:\ (Z,S)\in\chi\}$ \label{eq:sns_alm_subproblem2}
    \STATE $Y_{k+1}\gets -\grad f_\mu(X_{k+1})+\rho (X_{k+1}-Z_{k+1})$
    \STATE $k \gets k + 1$
    \ENDWHILE
    \end{algorithmic}
    }
\end{algorithm}

Using the fast version of \textbf{Algorithm~\ref{alg:alms}}, a fast version of \textbf{Algorithm~\ref{alg:sns_alm}} with $\cO(\rho/k^2)$ convergence rate, employing partial splitting and alternating linearization, can be constructed. This fast version can compute an $\epsilon$-optimal and feasible solution to problem $(P)$ in $\cO(1/\epsilon)$ iterations. Moreover, like the proximal gradient methods described earlier, each iteration for these methods can be computed efficiently. The subproblems to be solved at each iteration of \textbf{Algorithm~\ref{alg:sns_alm}} and its fast version have the following generic form:
\begin{align}
&\min_{X\in\reals^{m\times n}} f_\mu(X)+\fprod{Q, X-\tilde{X}}+\frac{\rho}{2}\norm{X-\tilde{X}}_F^2, \label{eq:sns_alm_subproblem1_generic}\\
&\min_{Z,S\in\reals^{m\times n}} \{\xi\norm{S}_1+\fprod{Q,Z-\tilde{Z}}+\frac{\rho}{2}\norm{Z-\tilde{Z}}_F^2:\ (Z,S)\in\chi\} \label{eq:sns_alm_subproblem2_generic}.
\end{align}
Let $U~\diag(\sigma)V^T$ denote the singular value decomposition of the matrix $\tilde{X}-Q/\rho$, then $X^*$, the minimizer of the subproblem in \eqref{eq:sns_alm_subproblem1_generic}, can be easily computed as $U~\diag\left(\sigma-\frac{\sigma}{\max\{\rho\sigma,\ 1+\rho\mu\}}\right)V^T$. And Lemma~\ref{lem:subproblem} shows how to solve the subproblem in \eqref{eq:sns_alm_subproblem2_generic}.
\section{Non-smooth Augmented Lagrangian Algorithm}
\label{sec:nsa}
\textbf{Algorithm~\ref{alg:nsa}} is a Non-Smooth Augmented Lagrangian Algorithm~(\proc{NSA}) that solves the non-smooth problem $(P)$.
The subproblem in Step~\ref{algeq:subproblem1} of \textbf{Algorithm~\ref{alg:nsa}} is a matrix shrinkage problem and can be solved efficiently by computing a singular value decomposition~(SVD) of an $m\times n$ matrix; and Lemma~\ref{lem:subproblem} shows that the subproblem in Step~\ref{algeq:subproblem2} can also be solved efficiently.

\begin{algorithm}[h!]
    \caption{NSA($Z_0,Y_0$)}\label{alg:nsa}
    {\small
    \begin{algorithmic}[1]
    \STATE \textbf{input:} $Z_0\in\reals^{m\times n}$, $Y_0\in\reals^{m\times n}$
    \STATE $k \gets 0$
    \WHILE{$k\leq 0$}
    \STATE $X_{k+1}\gets\argmin_X\{\norm{X}_*+\fprod{Y_k, X-Z_k}+\frac{\rho_k}{2}\norm{X-Z_k}_F^2\}$ \label{algeq:subproblem1}
    \STATE $\hat{Y}_{k+1}\gets Y_k+\rho_k (X_{k+1}-Z_{k})$
    \STATE $(Z_{k+1},S_{k+1})\gets\argmin_{\{(Z,S): \norm{Z+S-D}^2_F\leq\delta^2\}}\{\xi\norm{S}_1+\fprod{-Y_k, Z-X_{k+1}}+\frac{\rho_k}{2}\norm{Z-X_{k+1}}_F^2
    \}$ \label{algeq:subproblem2}
    \STATE Let $\theta_{k}$ be an optimal Lagrangian dual variable for the $\frac{1}{2}\norm{Z+S-D}^2_F\leq\frac{\delta^2}{2}$ constraint
    \STATE $Y_{k+1}\gets Y_k+\rho_k (X_{k+1}-Z_{k+1})$
    \STATE Choose $\rho_{k+1}$ such that $\rho_{k+1}\geq\rho_{k}$
    \STATE $k \gets k + 1$
    \ENDWHILE
    \end{algorithmic}
    }
\end{algorithm}

We now prove that Algorithm~\proc{NSA} converges under fairly mild conditions on the sequence $\{\rho_k\}_{k\in\integers_+}$ of penalty parameters. We first need the following lemma, which extends the similar result given in \cite{Ma09_1J} to partial splitting of variables.

\begin{lemma}
\label{lem:finite_sums}
Suppose that $\delta>0$. Let $\{X_k,Z_k,S_k,Y_k,\theta_k\}_{k\in\integers_+}$ be the sequence produced by Algorithm~\proc{NSA}. $(X^*,X^*,S^*)=\argmin_{X,Z,S}\{\norm{X}_*+\xi~\norm{S}_1:\ \frac{1}{2}\norm{Z+S-D}^2_F\leq\frac{\delta^2}{2},\ X=Z\}$ be any optimal solution, $Y^*\in\reals^{m\times n}$ and $\theta^*\geq 0$ be any optimal Lagrangian duals corresponding to the constraints $X=Z$ and $\frac{1}{2}\norm{Z+S-D}^2_F\leq\frac{\delta^2}{2}$, respectively. Then $\{\norm{Z_{k}-X^*}_F^2+\rho_{k}^{-2}\norm{Y_{k}-Y^*}_F^2\}_{k\in\integers_+}$ is a non-increasing sequence and
\begin{equation*}
\begin{array}{ll}
\sum_{k\in\integers_+}\norm{Z_{k+1}-Z_k}_F^2<\infty\hspace{5mm} &\sum_{k\in\integers_+}\rho_{k}^{-2}\norm{Y_{k+1}-Y_k}_F^2<\infty,\\
\sum_{k\in\integers_+}\rho_k^{-1}\fprod{-Y_{k+1}+Y^*, S_{k+1}-S^*}<\infty\hspace{5mm} &\sum_{k\in\integers_+}\rho_k^{-1}\fprod{-\hat{Y}_{k+1}+Y^*, X_{k+1}-X^*}<\infty,
\end{array}
\end{equation*}
\vspace{-5mm}
\begin{equation*}
\begin{array}{c}
\sum_{k\in\integers_+}\rho_k^{-1}\fprod{Y^*-Y_{k+1}, X^*+S^*-Z_{k+1}-S_{k+1}}<\infty.
\end{array}
\end{equation*}
\end{lemma}
\begin{proof}
See Appendix~\ref{sec:lem_proof} for the proof.
\end{proof}

Given partially split SPCP problem, $\min_{X,Z,S}\{\norm{X}_*+\xi\norm{S}_1:~X=Z,~(Z,S)\in\chi\}$, let $\mathcal{L}$ be its Lagrangian function
\begin{align}
\label{eq:lagrangian_split}
\mathcal{L}(X,Z,S;Y,\theta)=\norm{X}_*+\xi~\norm{S}_1+\fprod{Y,X-Z}+\frac{\theta}{2}\left(\norm{Z+S-D}_F^2-\delta^2\right).
\end{align}
\begin{theorem}
Suppose that $\delta>0$. Let $\{X_k,Z_k,S_k,Y_k,\theta_k\}_{k\in\integers_+}$ be the sequence produced by Algorithm~\proc{NSA}. Choose $\{\rho_k\}_{k\in\integers_+}$ such that
\begin{enumerate}[(i)]
\item $\sum_{k\in\integers_+}\frac{1}{\rho_k}=\infty$: Then $\lim_{k\in\integers_+}Z_k=\lim_{k\in\integers_+}X_k=X^*$, $\lim_{k\in\integers_+}S_k=S^*$ such that $(X^*, S^*)=\argmin\{\norm{X}_*+\xi~\norm{S}_1:\ \norm{X+S-D}_F\leq\delta\}$.
\item $\sum_{k\in\integers_+}\frac{1}{\rho_k^2}=\infty$: If $\norm{D-X^*}_F\neq \delta$, then $\lim_{k\in\integers_+}\theta_k=\theta^*\geq 0$ and $\lim_{k\in\integers_+}Y_k=Y^*$ such that $(X^*,X^*,S^*,Y^*,\theta^*)$ is a saddle point of the Lagrangian function $\mathcal{L}$ in \eqref{eq:lagrangian_split}. Otherwise, if $\norm{D-X^*}_F=\delta$, then there exists a limit point, $(Y^*,\theta^*)$, of the sequence $\{Y_k,\theta_k\}_{k\in\integers_+}$ such that $(Y^*,\theta^*)=\argmax_{Y,\theta}\{\mathcal{L}(X^*,X^*,S^*;Y,\theta):\ \theta\geq 0\}$.
\end{enumerate}
\end{theorem}
\begin{remark}
Requiring $\sum_{k\in\integers_+}\frac{1}{\rho_k}=\infty$ is similar to the condition in Theorem~2 in \cite{Ma09_1J}, which is needed to show that Algorithm I-ALM converges to an optimal solution of the robust PCA problem.
\end{remark}
\begin{remark}
Let $D=X^0+S^0+\zeta^0$ such that $\norm{\zeta^0}_F\leq\delta$ and $(X^0,S^0)$ satisfies the assumptions of Theorem~\ref{thm:candes2}. If $\norm{S^0}_F>\sqrt{Cmn}\delta$, then with very high probability, $\norm{D-X^*}_F>\delta$, where $C$ is the numerical constant defined in Theorem~\ref{thm:candes2}. Therefore, most of the time in applications, one does not encounter the case where $\norm{D-X^*}_F=\delta$.
\end{remark}
\begin{proof}
From Lemma~\ref{lem:finite_sums} and the fact that $X_{k+1}-Z_{k+1}=\frac{1}{\rho_k}~(Y_{k+1}-Y_k)$ for all $k\geq 1$, we have
\begin{align*}
\infty>\sum_{k\in\integers_+}\rho_{k}^{-2}\norm{Y_{k+1}-Y_k}_F^2=\sum_{k\in\integers_+}\norm{X_{k+1}-Z_{k+1}}_F^2.
\end{align*}
Hence, $\lim_{k\in\integers_+}(X_k-Z_k)=0$.

Let $(X^\#,X^\#,S^\#)=\argmin_{X,Z,S}\{\norm{X}_*+\xi~\norm{S}_1:\ \frac{1}{2}\norm{Z+S-D}^2_F\leq\frac{\delta^2}{2},\ X=Z\}$ be any optimal solution, $Y^\#\in\reals^{m\times n}$ and $\theta^\#\geq 0$ be any optimal Lagrangian duals corresponding to $X=Z$ and $\frac{1}{2}\norm{Z+S-D}^2_F\leq\frac{\delta^2}{2}$ constraints, respectively and $f^*:=\norm{X^\#}_*+\xi~\norm{S^\#}_1$.

Moreover, let $\chi=\{(Z,S)\in\reals^{m\times n}\times\reals^{m\times n}:\ \norm{Z+S-D}_F\leq\delta\}$ and $\mathbf{1}_\chi(Z,S)$ denote the indicator function of the closed convex set $\chi$, i.e. $\mathbf{1}_\chi(Z,S)=0$ if $(Z,S)\in\chi$; otherwise, $\mathbf{1}_\chi(Z,S)=\infty$. Since the sequence $\{(Z_k,S_k)\}_{k\in\integers_+}$ produced by \proc{NSA} is a feasible sequence for the set $\chi$, we have $\mathbf{1}_\chi(Z_k,S_k)=0$ for all $k\geq 1$. Hence, the following inequality is true for all $k\geq 0$
\begin{align}
 &\norm{X_k}_*+\xi~\norm{S_k}_1 \nonumber\\
=~&\norm{X_k}_*+\xi~\norm{S_k}_1 + \mathbf{1}_\chi(Z_k,S_k), \nonumber \\
\leq~ &\norm{X^\#}_*+\xi~\norm{S^\#}_1 + \mathbf{1}_\chi(X^\#,S^\#) -\fprod{-\hat{Y}_k, X^\#-X_k}-\fprod{-Y_k, S^\#-S_k}-\fprod{Y_k, X^\#+S^\#-Z_k-S_k}, \nonumber \\
=~& f^* + \fprod{-\hat{Y}_k+Y^\#, X_k-X^\#}+\fprod{-Y_k+Y^\#, S_k-S^\#}+\fprod{Y^\#-Y_k, X^\#+S^\#-Z_k-S_k} \label{eq:convexity_bound} \\
&+\fprod{Y^\#, Z_k-X_k},\nonumber
\end{align}
where the inequality follows from the convexity of norms and the fact that $-Y_k\in\xi~\partial\norm{S_k}_1$, $-\hat{Y}_k\in\partial\norm{X_k}_*$ and $(Y_k,Y_k)\in\partial\mathbf{1}_\chi(Z_k,S_k)$; the final equality follows from rearranging the terms and the fact that $(X^\#,S^\#)\in\chi$.

From Lemma~\ref{lem:finite_sums}, we have
$$\sum_{k\in\integers_+}\rho_{k-1}^{-1}\left(\fprod{-\hat{Y}_{k}+Y^\#, X_{k}-X^\#}+\fprod{-Y_{k}+Y^\#, S_{k}-S^\#}+\fprod{Y^\#-Y_{k}, X^\#+S^\#-Z_{k}-S_{k}}\right)<\infty.$$
Since $\sum_{k\in\integers_+}\frac{1}{\rho_k}=\infty$, there exists $\cK\subset\integers_+$ such that
\begin{align}
\lim_{k\in\cK}\left(\fprod{-\hat{Y}_{k}+Y^\#, X_{k}-X^\#}+\fprod{-Y_{k}+Y^\#, S_{k}-S^\#}+\fprod{Y^\#-Y_{k}, X^\#+S^\#-Z_{k}-S_{k}}\right)=0. \label{eq:inner_product_limit}
\end{align}
\eqref{eq:inner_product_limit} and the fact that $\lim_{k\in\integers_+}Z_k-X_k=0$ imply that along $\cK$ \eqref{eq:convexity_bound} converges to $f^*=\norm{X^\#}_*+\xi~\norm{S^\#}_1=\min\{\norm{X}_*+\xi~\norm{S}_1:\ (X,S)\in\chi\}$; hence along $\cK$ subsequence, $\{\norm{X_k}_*+\xi~\norm{S_k}_1\}_{k\in\cK}$ is a bounded sequence. Therefore, there exists $\cK^*\subset\cK\subset\integers_+$ such that $\lim_{k\in\cK^*}(X_k,S_k)=(X^*,S^*)$. Also, since $\lim_{k\in\integers_+}Z_k-X_k=0$ and $(Z_k,S_k)\in\chi$ for all $k\geq1$, we also have $(X^*,S^*)=\lim_{k\in\cK^*}(Z_k,S_k)\in\chi$. Since the limit of both sides of \eqref{eq:convexity_bound} along $\cK^*$ gives $\norm{X^*}_*+\xi~\norm{S^*}_1=\lim_{k\in\cK^*}\norm{X_k}_*+\xi~\norm{S_k}_1\leq f^*$ and $(X^*,S^*)\in\chi$, we conclude that $(X^*,S^*)=\argmin\{\norm{X}_*+\xi~\norm{S}_1:\ (X,S)\in\chi\}$.

It is also true that $(X^*,X^*,S^*)$ is an optimal solution to an equivalent problem: $\argmin_{X,Z,S}\{\norm{X}_*+\xi~\norm{S}_1:\ \frac{1}{2}\norm{Z+S-D}^2_F\leq\frac{\delta^2}{2},\ X=Z\}$. Now, let $\bar{Y}\in\reals^{m\times n}$ and $\bar{\theta}\geq 0$ be optimal Lagrangian duals corresponding to $X=Z$ and $\frac{1}{2}\norm{Z+S-D}^2_F\leq\frac{\delta^2}{2}$ constraints, respectively. From Lemma~\ref{lem:finite_sums}, it follows that $\{\norm{Z_{k}-X^*}_F^2+\rho_{k}^{-2}\norm{Y_{k}-\bar{Y}}_F^2\}_{k\in\integers_+}$ is a bounded non-increasing sequence. Hence, it has a unique limit point, i.e.
\begin{align*}
\lim_{k\in\integers_+}\norm{Z_{k}-X^*}_F^2=\lim_{k\in\integers_+}\norm{Z_{k}-X^*}_F^2+\rho_{k}^{-2}\norm{Y_{k}-\bar{Y}}_F^2=\lim_{k\in\cK^*}\norm{Z_{k}-X^*}_F^2+\rho_{k}^{-2}\norm{Y_{k}-\bar{Y}}_F^2=0,
\end{align*}
where the equalities follow from the facts that $\lim_{k\in\cK^*}Z_k=X^*$, $\mu_k\nearrow\infty$ as $k\rightarrow\infty$ and $\{\hat{Y}_k\}_{k\in\integers_+}$, $\{Y_k\}_{k\in\integers_+}$ are bounded sequences.
$\lim_{k\in\integers_+}\norm{Z_{k}-X^*}_F=0$ and $\lim_{k\in\integers_+} Z_{k}-X_k=0$ imply that $\lim_{k\in\integers_+}X_k=X^*$.

Using Lemma~\ref{lem:subproblem} for the $k$-th subproblem given in Step~\ref{algeq:subproblem2} in \textbf{Algorithm~\ref{alg:nsa}}, we have
\begin{align}
&S_{k+1}=sign\left(D-\left(X_{k+1}+\frac{1}{\rho_k}~Y_k\right)\right)\odot\max\left\{\left|D-\left(X_{k+1}+\frac{1}{\rho_k}~Y_k\right)\right|-\xi\frac{(\rho_k+\theta_k)}{\rho_k\theta_k}~E,\ \mathbf{0}\right\}, \label{eq:Sk}\\
&Z_{k+1}= \frac{\theta_k}{\rho_k+\theta_k}~(D-S_{k+1})+\frac{\rho_k}{\rho_k+\theta_k}~\left(X_{k+1}+\frac{1}{\rho_k}~Y_k\right). \label{eq:Xk}
\end{align}
If $\norm{D-(X_{k+1}+\frac{1}{\rho_k}~Y_k)}_F\leq\delta$, then $\theta_k=0$; otherwise, $\theta_k>0$ is the unique solution such that $\phi_k(\theta_k)=\delta$, where
\begin{align}
\phi_k(\theta):= \left\|\min\left\{\frac{\xi}{\theta}~E,\ \frac{\rho_k}{\rho_k+\theta}~\left|D-\left(X_{k+1}+\frac{1}{\rho_k}~Y_k\right)\right|\right\}\right\|_F.
\end{align}
In the following, it is shown that the sequence $\{S_k\}_{k\in\integers_+}$ has a unique limit point $S^*$. Since $\lim_{k\in\integers_+}X_k=X^*$, $\{Y_k\}_{k\in\integers_+}$ is a bounded sequence and $\rho_k\nearrow\infty$ as $k\rightarrow\infty$, we have $\lim_{k\in\integers_+}X_{k+1}+\frac{1}{\rho_k}~Y_k=X^*$.
\subsection*{Case 1: $\norm{D-X^*}_F\leq\delta$}
Previously, we have shown that that exists a subsequence $\cK^*\subset\integers_+$ such that $\lim_{k\in\cK^*}(X_k,S_k)=(X^*,S^*)=\argmin_{X,S}\{\norm{X}_*+\xi\norm{S}_1:\ \norm{X+S-D}_F\leq\delta\}$. On the other hand, since $\norm{D-X^*}_F\leq\delta$, $(X^*,\mathbf{0})$ is a feasible solution. Hence, $\norm{X^*}_*+\xi\norm{S^*}\leq\norm{X^*}_*$, which implies that $S^*=\mathbf{0}$.
\begin{align}
&\norm{X_k}_*+\xi~\norm{S_k}_1 \nonumber\\
=~&\norm{X_k}_*+\xi~\norm{S_k}_1 + \mathbf{1}_\chi(Z_k,S_k), \nonumber \\
\leq~&\norm{X^*}_*+\xi~\norm{\mathbf{0}}_1 + \mathbf{1}_\chi(X^*,\mathbf{0})-\fprod{-\hat{Y}_k, X^*-X_k}-\fprod{-Y_k, \mathbf{0}-S_k}-\fprod{Y_k, X^*+\mathbf{0}-Z_k-S_k},\nonumber\\
=~&\norm{X^*}_*+\fprod{\hat{Y}_k, X^*-X_k}+\fprod{Y_k, Z_k-X^*}. \label{eq:S_limit_equality_cond}
\end{align}
Since the sequences $\{Y_k\}_{k\in\integers_+}$ and $\{\hat{Y}_k\}_{k\in\integers_+}$ are bounded and $\lim_{k\in\integers_+}X_k=\lim_{k\in\integers_+}Z_k=X^*$, taking the limit on both sides of \eqref{eq:S_limit_equality_cond}, we have
\begin{align*}
  &\norm{X^*}_*+\xi~\lim_{k\in\integers_+}\norm{S_k}_1=\lim_{k\in\integers_+}\norm{X_k}_*+\xi~\norm{S_k}_1\\
=~&\lim_{k\in\integers_+}\norm{X^*}_*+\fprod{\hat{Y}_k, X^*-X_k}+\fprod{Y_k, Z_k-X^*}=\norm{X^*}_*.
\end{align*}
Therefore, $\lim_{k\in\integers_+}\norm{S_k}_1=0$, which implies that $\lim_{k\in\integers_+}S_k=S^*=\mathbf{0}$.

\subsection*{Case 2: $\norm{D-X^*}_F>\delta$}
Since $\norm{D-(X_{k+1}+\frac{1}{\rho_k}~Y_k)}_F\rightarrow \norm{D-X^*}_F>\delta$, there exists $K\in\integers_+$ such that for all $k\geq K$, $\norm{D-(X_{k+1}+\frac{1}{\rho_k}~Y_k)}_F>\delta$. For all $k\geq K$, $\phi_k(.)$ is a continuous and strictly decreasing function of $\theta$ for $\theta\geq 0$. Hence, inverse function $\phi^{-1}_k(.)$ exits around $\delta$ for all $k\geq K$. Thus, $\phi_k(0)=\norm{D-(X_{k+1}+\frac{1}{\rho_k}~Y_k)}_F>\delta$ and $\lim_{\theta\rightarrow\infty}\phi_k(\theta)=0$ imply that $\theta_k=\phi^{-1}_k(\delta)>0$ for all $k\geq K$. Moreover, $\phi_k(\theta)\leq\phi(\theta):=\norm{\frac{\xi}{\theta}~E}_F$ implies that $\theta_k\leq\frac{\xi\sqrt{mn}}{\delta}$ for all $k\geq K$. Therefore, $\{\theta_k\}_{k\in\integers_+}$ is a bounded sequence, which has a convergent subsequence $\cK_\theta\subset\integers_+$ such that $\lim_{k\in\cK_\theta}\theta_k=\theta^*$. We also have $\phi_k(\theta)\rightarrow\phi_\infty(\theta)$ pointwise for all $0\leq\theta\leq\frac{\xi\sqrt{mn}}{\delta}$, where
\begin{align}
\phi_\infty(\theta):= \left\|\min\left\{\frac{\xi}{\theta}~E,\ \left|D-X^*\right|\right\}\right\|_F.
\end{align}
Since $\phi_k(\theta_k)=\delta$ for all $k\geq K$, we have
\begin{align}
\delta=\lim_{k\in\cK}\phi_k(\theta_k)=\left\|\min\left\{\frac{\xi}{\theta_k}~E,\ \frac{\rho_k}{\rho_k+\theta_k}~\left|D-\left(X_{k+1}+\frac{1}{\rho_k}~Y_k\right)\right|\right\}\right\|_F=\phi_\infty(\theta^*).
\end{align}
Note that since $\norm{D-X^*}_F>\delta$, $\phi_\infty$ is invertible around $\delta$, i.e. $\phi_\infty^{-1}$ exists around $\delta$. Thus, $\theta^*=\phi_\infty^{-1}(\delta)$. Since $\cK_\theta$ is an arbitrary subsequence, we can conclude that $\theta^*:=\lim_{k\in\integers_+}\theta_k=\phi_\infty^{-1}(\delta)$. Since there exists $\theta^*>0$ such that $\theta^*=\lim_{k\in\integers_+}\theta_k$, taking the limit on both sides of \eqref{eq:Sk}, we have
\begin{align}
S^*:=\lim_{k\in\integers_+}S_{k+1}=sign\left(D-X^*\right)\odot\max\left\{\left|D-X^*\right|-\frac{\xi}{\theta^*}~E,\ \mathbf{0}\right\},
\end{align}
and this completes the first part of the theorem.

Now, we will show that if $\norm{D-X^*}_F\neq\delta$, then the sequences $\{\theta_k\}_{k\in\integers_+}$ and $\{Y_k\}_{k\in\integers_+}$ have unique limits. Note that from \eqref{eq:Zopt_cond}, it follows that $Y_k=\theta_{k-1}(Z_k+S_k-D)$ for all $k\geq 1$. First suppose that $\norm{D-X^*}_F<\delta$. Since $\norm{D-(X_{k+1}+\frac{1}{\rho_k}~Y_k)}_F\rightarrow \norm{D-X^*}_F<\delta$, there exists $K\in\integers_+$ such that for all $k\geq K$, $\norm{D-(X_{k+1}+\frac{1}{\rho_k}~Y_k)}_F<\delta$. Thus, from Lemma~\ref{lem:subproblem} for all $k\geq K$, $\theta_k=0$, $S_{k+1}=0$, $Z_{k+1}=X_{k+1}+\frac{1}{\rho_k}~Y_k$, which implies that $\theta^*:=\lim_{k\in\integers_+}\theta_k=0$ and $Y^*=\lim_{k\in\integers_+}Y_k=\lim_{k\in\integers_+}\theta_{k-1}(Z_{k}+S_{k}-D)=\mathbf{0}$ since $S^*=\lim_{k\in\cK^*}S_k=\lim_{k\in\integers_+}S_k=0$, $\lim_{k\in\integers_+}Z_k=X^*$ and $\norm{D-X^*}_F<\delta.$ Now suppose that $\norm{D-X^*}_F>\delta$. In Case~2 above we have shown that $\theta^*=\lim_{k\in\integers_+}\theta_k$. Hence, there exists $Y^*\in\reals^{m\times n}$ such that $Y^*=\lim_{k\in\integers_+}\theta_{k-1}(Z_k+S_k-D)=\theta^*(X^*+S^*-D)$.


Suppose that $\sum_{k\in\integers_+}\frac{1}{\rho_k^2}=\infty$. From Lemma~\ref{lem:finite_sums}, we have $\sum_{k\in\integers_+}\norm{Z_{k+1}-Z_k}_F^2<\infty$. Equivalently, the series can be written as
\begin{align}
\infty>\sum_{k\in\integers_+}\norm{Z_{k+1}-Z_k}_F^2=\sum_{k\in\integers_+}\rho_k^{-2}\norm{\hat{Y}_{k+1}-Y_{k+1}}_F^2.
\end{align}
Since $\sum_{k\in\integers_+}\frac{1}{\rho_k^2}=\infty$, there exists a subsequence $\cK\subset\integers_+$ such that $\lim_{k\in\cK}\norm{\hat{Y}_{k+1}-Y_{k+1}}_F^2=0$. Hence, $\lim_{k\in\cK}\rho_k^2\norm{Z_{k+1}-Z_k}_F^2=0$, i.e. $\lim_{k\in\cK}\rho_k(Z_{k+1}-Z_k)=0$.

Using \eqref{eq:Xopt_cond}, \eqref{eq:Sopt_cond} and \eqref{eq:Zopt_cond}, we have
\begin{align}
&0\in\partial\norm{X_{k+1}}_*+\theta_k(Z_{k+1}+S_{k+1}-D)+\rho_k(Z_{k+1}-Z_k), \label{eq:Xk_opt}\\
&0\in\xi\partial\norm{S_{k+1}}_1+ \theta_k(Z_{k+1}+S_{k+1}-D).\label{eq:Sk_opt}
\end{align}

If $\norm{D-X^*}\neq\delta$,  then there exists $Y^*\in\reals^{m\times n}$ such that $Y^*=\lim_{k\in\integers_+}\theta_{k-1}(Z_k+S_k-D)=\theta^*(X^*+S^*-D)$.
Taking the limit of \eqref{eq:Xk_opt},\eqref{eq:Sk_opt} along $\cK\subset\integers_+$ and using the fact that $\lim_{k\in\cK}\rho_k(Z_{k+1}-Z_k)=0$, we have
\begin{align}
&0\in\partial\norm{X^*}_*+\theta^*(X^*+S^*-D), \label{eq:X_opt}\\
&0\in\xi\partial\norm{S^*}_1+ \theta^*(X^*+S^*-D).\label{eq:S_opt}
\end{align}
\eqref{eq:X_opt} and \eqref{eq:S_opt} together imply that $(X^*, S^*)$, $Y^*=\theta^*(X^*+S^*-D)$ and $\theta^*$ satisfy KKT optimality conditions for the problem $\min_{X,Z,S}\{\norm{X}_*+\xi~\norm{S}_1:\ \frac{1}{2}\norm{Z+S-D}^2_F\leq\frac{\delta^2}{2},\ X=Z\}$. Hence, $(X^*,X^*,S^*,Y^*,\theta^*)$ is a saddle point of the Lagrangian function
\begin{align*}
\mathcal{L}(X,Z,S;Y,\theta)=\norm{X}_*+\xi~\norm{S}_1+\fprod{Y,X-Z}+\frac{\theta}{2}\left(\norm{Z+S-D}_F^2-\delta^2\right).
\end{align*}

Suppose that $\norm{D-X^*}_F=\delta$. Fix $k>0$. If $\norm{D-(X_{k+1}+\frac{1}{\rho_k}~Y_k)}_F\leq \delta$, then $\theta_k=0$. Otherwise, $\theta_k>0$ and as shown in case~2 in the first part of the proof $\theta_k\leq\frac{\xi\sqrt{mn}}{\delta}$. Thus, for any $k>0$, $0\leq\theta_k\leq\frac{\xi\sqrt{mn}}{\delta}$. Since $\{\theta_k\}_{k\in\integers_+}$ is a bounded sequence, there exists a further subsequence $\cK_\theta\subset\cK$ such that $\theta^*:=\lim_{k\in\cK_\theta}\theta_{k-1}$ and $Y^*:=\lim_{k\in\cK_\theta}\theta_{k-1}(Z_k+S_k-D)=\theta^*(X^*+S^*-D)$ exist. Thus, taking the limit of \eqref{eq:Xk_opt},\eqref{eq:Sk_opt} along $\cK_\theta\subset\integers_+$ and using the facts that $\lim_{k\in\cK}\rho_k(Z_{k+1}-Z_k)=0$ and $X^*=\lim_{k\in\integers_+}X_k=\lim_{k\in\integers_+}Z_k$, $S^*=\lim_{k\in\integers_+}S_k$ exist, we conclude that $(X^*,X^*,S^*,Y^*,\theta^*)$ is a saddle point of the Lagrangian function $\mathcal{L}(X,Z,S;Y,\theta)$.
\end{proof}
\section{Numerical experiments}
\label{sec:computations}
Our preliminary numerical experiments showed that among the four algorithms discussed in this paper, \proc{NSA} is the fastest. It also has very few parameters that need to be tuned. Therefore, 
we only report the results for \proc{NSA}. We conducted two sets of numerical experiments with \proc{NSA} to solve ~\eqref{eq:stable_component_pursuit}, where $\xi=\frac{1}{\sqrt{\max\{m,n\}}}$. In the first set we solved randomly generated instances of the stable principle component pursuit problem. In this setting, first we tested only \proc{NSA} to see how the run times scale with respect to problem parameters and size; then we compared \proc{NSA} with another alternating direction augmented Lagrangian algorithm ASALM~\cite{Tao09_1J}. In the second set of experiments, we ran NSA and ASALM to extract moving objects from an airport security noisy video~\cite{Li04_1J}.
\subsection{Random Stable Principle Component Pursuit Problems}
\label{sec:rPCA_results}
We tested \proc{NSA} on randomly generated stable principle component pursuit problems. The data matrices for these problems, $D=X^0+S^0+\zeta^0$, were generated as follows
\begin{enumerate}[i.]
\item $X^0=UV^T$, such that $U\in\reals^{n\times r}$, $V\in\reals^{n\times r}$ for $r=c_r n$ and $U_{ij}\sim \mathcal{N}(0,1)$, $V_{ij}\sim
  \mathcal{N}(0,1)$ for all $i,j$ are independent standard Gaussian variables and $c_r\in\{0.05, 0.1\}$,
\item $\Lambda\subset\{(i,j):\ 1 \leq i,j\leq n\}$ such that cardinality of $\Lambda$, $|\Lambda|=p$ for $p=c_p n^2$ and $c_p\in\{0.05, 0.1\}$,
\item $S^0_{ij}\sim\mathcal{U}[-100,100]$ for all $(i,j)\in\Lambda$ are independent uniform random variables between $-100$ and $100$,
\item $\zeta^0_{ij}\sim \varrho\mathcal{N}(0,1)$ for all $i,j$ are independent Gaussian variables.
\end{enumerate}
We created 10 random problems of size $n\in\{500, 1000, 1500\}$, i.e. $D\in\reals^{n
  \times n}$, for each of the two choices of $c_r$ and $c_p$ using the procedure described above, where $\varrho$ was set such that signal-to-noise ratio of $D$ is either $80dB$ or $45dB$. Signal-to-noise ratio of $D$ is given by
  \begin{align}
  \label{eq:snr}
  \rm{SNR}(D)=10\log_{10}\left(\frac{E\left[\norm{X^0+S^0}_F^2\right]}{E\left[\norm{\zeta^0}_F^2\right]}\right)=10\log_{10}\left(\frac{c_r n+c_s 100^2/3}{\varrho^2}\right).
  \end{align}
  Hence, for a given SNR value, we selected $\varrho$ according to \eqref{eq:snr}. Table~\ref{tab:snr} displays the $\varrho$ value we have used in our experiments.
  \begin{table}[!htb]
    \begin{adjustwidth}{-2em}{-2em}
    \centering
    \caption{$\varrho$ values depending on the experimental setting}
    \renewcommand{\arraystretch}{1.1}
    {\footnotesize
    \begin{tabular}{|c|c|c|c|c|c|c|}
    \hline
    SNR&n&$\mathbf{c_r}$=\textbf{0.05} $\mathbf{c_p}$=\textbf{0.05}&$\mathbf{c_r}$=\textbf{0.05} $\mathbf{c_p}$=\textbf{0.1}&$\mathbf{c_r}$=\textbf{0.1} $\mathbf{c_p}$=\textbf{0.05}&$\mathbf{c_r}$=\textbf{0.1} $\mathbf{c_p}$=\textbf{0.1}\\\hline
    \multirow{3}{*}{$80dB$}
    &$\mathbf{500}$
    & 0.0014 & 0.0019 & 0.0015 & 0.0020 \\ \cline{2-6}
    &$\mathbf{1000}$
    & 0.0015 & 0.0020 & 0.0016 & 0.0021 \\ \cline{2-6}
    &$\mathbf{1500}$
    & 0.0016 & 0.0020 & 0.0018 & 0.0022 \\ \hline
    \multirow{3}{*}{$45dB$}
    &$\mathbf{500}$
    & 0.0779 & 0.1064 & 0.0828 & 0.1101 \\ \cline{2-6}
    &$\mathbf{1000}$
    & 0.0828 & 0.1101 & 0.0918 & 0.1171 \\ \cline{2-6}
    &$\mathbf{1500}$
    & 0.0874 & 0.1136 & 0.1001 & 0.1236 \\ \hline
    \end{tabular}
    \label{tab:snr}
    }
    \end{adjustwidth}
\end{table}
As in \cite{Tao09_1J}, we set $\delta = \sqrt{(n + \sqrt{8n})}\varrho$ in \eqref{eq:stable_component_pursuit} in the first set of experiments for both \proc{NSA} and \proc{ASALM}.

Our code for \proc{NSA} was written in MATLAB 7.2 and can be found at~\url{http://www.columbia.edu/~nsa2106}. 
We terminated the algorithm when
\begin{align}
\label{eq:stopping_cond}
\frac{\norm{(X_{k+1},S_{k+1})-(X_{k},S_{k})}_F}{\norm{(X_{k},S_{k})}_F+1}\leq \varrho.
\end{align}
The results of our experiments are displayed in Tables~\ref{tab:self_time} and \ref{tab:self_quality}.
In Table~\ref{tab:self_time}, the row labeled $\mathbf{CPU}$ lists the running time of \proc{NSA} in
\emph{seconds} and the row labeled $\mathbf{SVD\#}$ lists the number of partial singular value decomposition~(SVD) computed by \proc{NSA}.
The minimum, average and maximum CPU times and number of partial SVD taken over the $10$ random
instances are given for each choice of $n$, $c_r$ and $c_p$ values. Table~\ref{tab:self_detail_80dB} and Table~\ref{tab:self_detail_45dB} in the appendix list additional error statistics.

With the stopping condition given in \eqref{eq:stopping_cond}, the solutions produced by \proc{NSA} have $\frac{\norm{X^{sol}+S^{sol}-D}_F}{\norm{D}_F}$ approximately $1.5\times 10^{-4}$ when ${\rm SNR}(D)=80dB$ and $5\times 10^{-3}$ when ${\rm SNR}(D)=45dB$, regardless of the problem dimension $n$ and the problem parameters related to the rank and sparsity of $D$, i.e. $c_r$ and $c_p$. After thresholding the singular values of $X^{sol}$ that were less than $1\times 10^{-12}$, \proc{NSA} found the true rank in all 120 random problems solved when ${\rm SNR}(D)=80dB$, and it found the true rank for 113 out of 120 problems when ${\rm SNR}(D)=45dB$, while for 6 of the remaining problems $\rank(X^{sol})$ is off from $\rank(X^0)$ only by 1. Table~\ref{tab:self_time} shows that the number of partial SVD was a very slightly increasing function of $n$, $c_r$ and $c_p$. Moreover, Table~\ref{tab:self_quality} shows that the relative error of the solution $(X^{sol},S^{sol})$ was almost constant for different $n$, $c_r$ and $c_p$ values.
\begin{table}[!htb]
    \begin{adjustwidth}{-2em}{-2em}
    \centering
    \caption{NSA: Solution time for decomposing $D\in\reals^{n\times n}$, $n\in\{500, 1000, 1500\}$}
    \renewcommand{\arraystretch}{1.3}
    {\footnotesize
    \begin{tabular}{ccc|c|c|c|c|}
    \cline{4-7}
    &&&$\mathbf{c_r}$=\textbf{0.05} $\mathbf{c_p}$=\textbf{0.05}&$\mathbf{c_r}$=\textbf{0.05} $\mathbf{c_p}$=\textbf{0.1}&$\mathbf{c_r}$=\textbf{0.1} $\mathbf{c_p}$=\textbf{0.05}&$\mathbf{c_r}$=\textbf{0.1} $\mathbf{c_p}$=\textbf{0.1}\\\hline
    \multicolumn{1}{|c|}{SNR}&\multicolumn{1}{|c|}{n}&Field& min/\textbf{avg}/max & min/\textbf{avg}/max & min/\textbf{avg}/max & min/\textbf{avg}/max\\ \hline
    \multicolumn{1}{|c|}{\multirow{6}{*}{$80dB$}}
    &\multicolumn{1}{|c|}{\multirow{2}{*}{$\mathbf{500}$}}
    & $\mathbf{SVD\#}$ & 9/\textbf{9.0}/9 & 9/\textbf{9.5}/10 & 10/\textbf{10.0}/10 & 11/\textbf{11}/11 \\
    \multicolumn{1}{|c|}{}&\multicolumn{1}{|c|}{}& $\mathbf{CPU}$ & 3.2/\textbf{4.4}/5.1 & 3.6/\textbf{5.1}/6.6 & 4.3/\textbf{5.2}/6.4& 5.0/\textbf{6.2}/8.1 \\ \cline{2-7}
    \multicolumn{1}{|c|}{}
    &\multicolumn{1}{|c|}{\multirow{2}{*}{$\mathbf{1000}$}}
    & $\mathbf{SVD\#}$ & 9/\textbf{9.9}/10 & 10/\textbf{10.0}/10 & 11/\textbf{11}/11& 12/\textbf{12.0}/12 \\
    \multicolumn{1}{|c|}{}&\multicolumn{1}{|c|}{}& $\mathbf{CPU}$ & 16.5/\textbf{19.6}/22.4 & 14.6/\textbf{20.7}/24.3 & 25.2/\textbf{26.9}/29.1	& 27.9/\textbf{31.2}/36.3 \\ \cline{2-7}
    \multicolumn{1}{|c|}{}
    &\multicolumn{1}{|c|}{\multirow{2}{*}{$\mathbf{1500}$}}
    & $\mathbf{SVD\#}$ & 10/\textbf{10.0}/10	& 10/\textbf{10.9}/11 & 12/\textbf{12.0}/12 & 12/\textbf{12.2}/13 \\
    \multicolumn{1}{|c|}{}&\multicolumn{1}{|c|}{}& $\mathbf{CPU}$ & 38.6/\textbf{44.1}/46.6 & 43.7/\textbf{48.6}/51.9 & 78.6/\textbf{84.1}/90.8 & 80.7/\textbf{97.7}/155.2 \\ \hline
    \multicolumn{1}{|c|}{\multirow{6}{*}{$45dB$}}
    &\multicolumn{1}{|c|}{\multirow{2}{*}{$\mathbf{500}$}}
    & $\mathbf{SVD\#}$ & 6/\textbf{6}/6 & 6/\textbf{6.9}/7 & 7/\textbf{7.1}/8 & 8/\textbf{8}/8 \\
    \multicolumn{1}{|c|}{}&\multicolumn{1}{|c|}{}& $\mathbf{CPU}$ & 2.3/\textbf{2.9}/4.2 & 2.9/\textbf{3.6}/4.5 & 2.9/\textbf{3.9}/6.2& 3.5/\textbf{4.2}/6.0 \\ \cline{2-7}
    \multicolumn{1}{|c|}{}
    &\multicolumn{1}{|c|}{\multirow{2}{*}{$\mathbf{1000}$}}
    & $\mathbf{SVD\#}$ & 7/\textbf{7.0}/7 & 7/\textbf{7.0}/7 & 8/\textbf{8.1}/9& 9/\textbf{9.0}/9 \\
    \multicolumn{1}{|c|}{}&\multicolumn{1}{|c|}{}& $\mathbf{CPU}$ & 11.5/\textbf{13.4}/17.4 & 10.6/\textbf{13.3}/17.9 & 17.1/\textbf{18.7}/20.7	& 19.7/\textbf{23.8}/28.9 \\ \cline{2-7}
    \multicolumn{1}{|c|}{}
    &\multicolumn{1}{|c|}{\multirow{2}{*}{$\mathbf{1500}$}}
    & $\mathbf{SVD\#}$ & 7/\textbf{7.9}/8	& 8/\textbf{8.0}/8 & 9/\textbf{9.0}/9 & 9/\textbf{9.0}/9 \\
    \multicolumn{1}{|c|}{}&\multicolumn{1}{|c|}{}& $\mathbf{CPU}$ & 34.1/\textbf{37.7}/44.0 & 30.7/\textbf{37.1}/45.6 & 55.6/\textbf{59.0}/63.7 & 55.9/\textbf{59.7}/64.8 \\ \hline
    \end{tabular}
    \label{tab:self_time}
    }
    \end{adjustwidth}
\end{table}
\begin{table}[!htb]
    \begin{adjustwidth}{-2em}{-2em}
    \centering
    \caption{NSA: Solution accuracy for decomposing $D\in\reals^{n\times n}$, $n\in\{500, 1000, 1500\}$}
    \renewcommand{\arraystretch}{1.5}
    {\footnotesize
    \begin{tabular}{ccc|c|c|c|c|}
    \cline{4-7}
    &&&$\mathbf{c_r}$=\textbf{0.05} $\mathbf{c_p}$=\textbf{0.05}&$\mathbf{c_r}$=\textbf{0.05} $\mathbf{c_p}$=\textbf{0.1}&$\mathbf{c_r}$=\textbf{0.1} $\mathbf{c_p}$=\textbf{0.05}&$\mathbf{c_r}$=\textbf{0.1} $\mathbf{c_p}$=\textbf{0.1}\\\hline
    \multicolumn{1}{|c|}{SNR}&\multicolumn{1}{|c|}{n}&Relative Error& \textbf{avg}~/~max & \textbf{avg}~/~max & \textbf{avg}~/~max & \textbf{avg}~/~max\\ \hline
    \multicolumn{1}{|c|}{\multirow{6}{*}{$80dB$}}
    &\multicolumn{1}{|c|}{\multirow{2}{*}{$\mathbf{500}$}}
    & $\mathbf{\frac{\norm{X^{sol}-X^0}_F}{\norm{X^0}_F}}$ & \textbf{4.0E-4}~/~4.2E-4 & \textbf{5.8E-4}~/~8.5E-4 & \textbf{3.6E-4}~/~3.9E-4 & \textbf{4.4E-4}~/~4.5E-4 \\
    \multicolumn{1}{|c|}{}&\multicolumn{1}{|c|}{}& $\mathbf{\frac{\norm{S^{sol}-S^0}_F}{\norm{S^0}_F}}$ & \textbf{1.7E-4}~/~1.8E-4 & \textbf{1.6E-4}~/~2.5E-4 & \textbf{1.6E-4}~/~1.8E-4 & \textbf{1.3E-4}~/~1.3E-4 \\ \cline{2-7}
    \multicolumn{1}{|c|}{}
    &\multicolumn{1}{|c|}{\multirow{2}{*}{$\mathbf{1000}$}}
    & $\mathbf{\frac{\norm{X^{sol}-X^0}_F}{\norm{X^0}_F}}$ & \textbf{2.0E-4}~/~2.4E-4 & \textbf{3.8E-4}~/~4.1E-4 & \textbf{2.2E-4}~/~2.2E-4 & \textbf{2.8E-4}~/~2.9E-4 \\
    \multicolumn{1}{|c|}{}&\multicolumn{1}{|c|}{}& $\mathbf{\frac{\norm{S^{sol}-S^0}_F}{\norm{S^0}_F}}$ & \textbf{1.2E-4}~/~1.4E-4 & \textbf{1.5E-4}~/~1.6E-4 & \textbf{1.2E-4}~/~1.3E-4 & \textbf{1.1E-4}~/~1.1E-4 \\ \cline{2-7}
    \multicolumn{1}{|c|}{}
    &\multicolumn{1}{|c|}{\multirow{2}{*}{$\mathbf{1500}$}}
    & $\mathbf{\frac{\norm{X^{sol}-X^0}_F}{\norm{X^0}_F}}$ & \textbf{1.8E-4}~/~2.2E-4 & \textbf{2.1E-4}~/~2.6E-4 & \textbf{1.3E-4}~/~1.3E-4 & \textbf{2.8E-4}~/~2.9E-4 \\
    \multicolumn{1}{|c|}{}&\multicolumn{1}{|c|}{}& $\mathbf{\frac{\norm{S^{sol}-S^0}_F}{\norm{S^0}_F}}$ & \textbf{1.3E-4}~/~1.6E-4 & \textbf{9.6E-5}~/~1.1E-4 & \textbf{8.1E-5}~/~8.5E-5 & \textbf{1.3E-4}~/~1.4E-4 \\ \hline
    \multicolumn{1}{|c|}{\multirow{6}{*}{$45dB$}}
    &\multicolumn{1}{|c|}{\multirow{2}{*}{$\mathbf{500}$}}
    & $\mathbf{\frac{\norm{X^{sol}-X^0}_F}{\norm{X^0}_F}}$ & \textbf{6.0E-3}~/~6.2E-3 & \textbf{8.0E-3}~/~9.2E-3 & \textbf{6.1E-3}~/~6.3E-3 & \textbf{8.1E-3}~/~8.2E-3 \\
    \multicolumn{1}{|c|}{}&\multicolumn{1}{|c|}{}& $\mathbf{\frac{\norm{S^{sol}-S^0}_F}{\norm{S^0}_F}}$ & \textbf{2.1E-3}~/~2.2E-3 & \textbf{2.3E-3}~/~2.7E-3 & \textbf{2.2E-3}~/~2.3E-3 & \textbf{2.7E-3}~/~2.9E-3 \\ \cline{2-7}
    \multicolumn{1}{|c|}{}
    &\multicolumn{1}{|c|}{\multirow{2}{*}{$\mathbf{1000}$}}
    & $\mathbf{\frac{\norm{X^{sol}-X^0}_F}{\norm{X^0}_F}}$ & \textbf{4.1E-3}~/~4.2E-3 & \textbf{6.1E-3}~/~6.2E-3 & \textbf{4.6E-3}~/~4.7E-3 & \textbf{6.0E-3}~/~6.5E-3 \\
    \multicolumn{1}{|c|}{}&\multicolumn{1}{|c|}{}& $\mathbf{\frac{\norm{S^{sol}-S^0}_F}{\norm{S^0}_F}}$ & \textbf{1.9E-3}~/~1.9E-3 & \textbf{2.4E-3}~/~2.5E-3 & \textbf{2.3E-3}~/~3.5E-3 & \textbf{3.1E-3}~/~3.7E-3 \\ \cline{2-7}
    \multicolumn{1}{|c|}{}
    &\multicolumn{1}{|c|}{\multirow{2}{*}{$\mathbf{1500}$}}
    & $\mathbf{\frac{\norm{X^{sol}-X^0}_F}{\norm{X^0}_F}}$ & \textbf{3.4E-3}~/~3.6E-3 & \textbf{4.7E-3}~/~4.7E-3 & \textbf{3.9E-3}~/~4.0E-3 & \textbf{5.3E-3}~/~5.3E-3 \\
    \multicolumn{1}{|c|}{}&\multicolumn{1}{|c|}{}& $\mathbf{\frac{\norm{S^{sol}-S^0}_F}{\norm{S^0}_F}}$ & \textbf{1.8E-3}~/~1.8E-3 & \textbf{2.3E-3}~/~2.3E-3 & \textbf{2.6E-3}~/~3.5E-3 & \textbf{3.1E-3}~/~3.1E-3 \\ \hline
    \end{tabular}
    \label{tab:self_quality}
    }
    \end{adjustwidth}
\end{table}

Next, we compared NSA with ASALM~\cite{Tao09_1J} for a fixed problem size, i.e. $n=1500$ where $D\in\reals^{n\times n}$. In all the numerical experiments, we terminated \proc{NSA} according to \eqref{eq:stopping_cond}. For random problems with ${\rm SNR}(D)=80dB$, we terminated \proc{ASALM} according to \eqref{eq:stopping_cond}. However, for random problems with ${\rm SNR}(D)=45dB$,  \proc{ASALM} produced solutions with $99\%$ relative errors when \eqref{eq:stopping_cond} was used. Therefore, for random problems with ${\rm SNR}(D)=45dB$, we terminated \proc{ASALM} either when it computed a solution with better relative errors comparing to \proc{NSA} solution for the same problem or when an iterate satisfied \eqref{eq:stopping_cond} with the righthand side replaced by $0.1\varrho$. The code for \proc{ASALM} was obtained from the authors of~\cite{Tao09_1J}.

The comparison results are displayed in Table~\ref{tab:compare_time} and Table~\ref{tab:compare_quality}. In Table~\ref{tab:compare_time}, the row labeled $\mathbf{CPU}$ lists the running time of each algorithm in \emph{seconds} and the row labeled $\mathbf{SVD\#}$ lists the number of partial SVD computation of each algorithm. In Table~\ref{tab:compare_time},
the minimum, average and maximum of CPU times and the number of partial SVD computation of each algorithm taken over the $10$ random instances are given for each two choices of $c_r$ and $c_p$. Moreover, Table~\ref{tab:compare_detail_80dB} and Table~\ref{tab:compare_detail_45dB} given in the appendix list different error statistics.

We used PROPACK~\cite{propack} for computing partial singular value decompositions. In order to estimate the rank of $X^0$, we followed the scheme proposed in Equation~(17) in \cite{Ma09_1J}.

\begin{figure} [!h]
\centering
\caption{NSA: Comparison of randomly selected 1500 components of $\zeta^0$ with absolute errors of those components in $S^{sol}$ and $\sigma(X^{sol})$. $D\in\reals^{n\times n}$, $n=1500$, ${\rm SNR}(D)=45dB$}
\hspace{-0.5cm}
\includegraphics[scale=0.38]{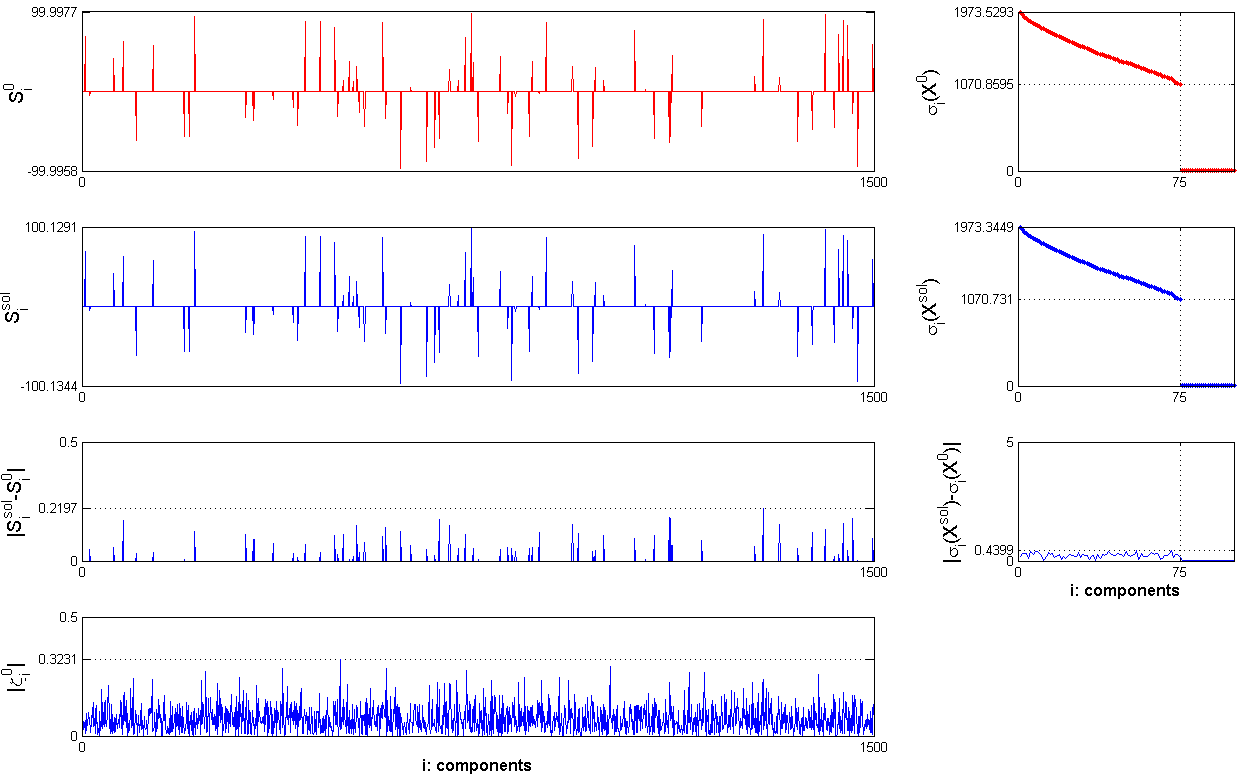}
\label{fig:nsa_45dB}
\end{figure}

Both \proc{NSA} and \proc{ASALM} found the true rank in all 40 random problems solved when ${\rm SNR}(D)=80dB$. \proc{NSA} found the true rank for 39 out of 40 problems with $n=1500$ when ${\rm SNR}(D)=45dB$, while for the remaining 1 problem $\rank(X^{sol})$ is off from $\rank(X^0)$ only by 1. On the other hand, when ${\rm SNR}(D)=45dB$, \proc{ASALM} could not find the true rank in any of the test problems. For each of the four problem settings corresponding to different $c_r$ and $c_p$ values, in Table~\ref{tab:rank} we report the average and maximum of $\rank(X^{sol})$ over 10 random instances, after thresholding the singular values of $X^{sol}$ that were less than $1\times 10^{-12}$.
\begin{table}[!htb]
    \begin{adjustwidth}{-2em}{-2em}
    \centering
    \caption{NSA vs ASALM: $\rank(X^{sol})$ values for problems with $n=1500$, ${\rm SNR}(D)=45dB$}
    \renewcommand{\arraystretch}{1.3}
    {\footnotesize
    \begin{tabular}{c|c|c|c|c|}
    \cline{2-5}
    &\multicolumn{2}{|c|}{$\rank(X^0)=75$}&\multicolumn{2}{|c|}{$\rank(X^0)=150$}\\ \cline{2-5}
    &$\mathbf{c_r}$=\textbf{0.05} $\mathbf{c_p}$=\textbf{0.05}&$\mathbf{c_r}$=\textbf{0.05} $\mathbf{c_p}$=\textbf{0.1}&$\mathbf{c_r}$=\textbf{0.1} $\mathbf{c_p}$=\textbf{0.05}&$\mathbf{c_r}$=\textbf{0.1} $\mathbf{c_p}$=\textbf{0.1}\\\hline
    \multicolumn{1}{|c|}{Alg.}& \textbf{avg}~/~max & \textbf{avg}~/~max & \textbf{avg}~/~max & \textbf{avg}~/~max\\ \hline
    \multicolumn{1}{|c|}{$\mathbf{NSA}$}
    &\textbf{75}~/~75&\textbf{75}~/~75&\textbf{150.1}~/~151&\textbf{150}~/~150\\ \hline
    \multicolumn{1}{|c|}{$\mathbf{ASALM}$}
    &\textbf{175.8}~/~177&\textbf{179}~/~207&\textbf{222.4}~/~224&\textbf{201.9}~/~204\\ \hline
    \end{tabular}
    \label{tab:rank}
    }
    \end{adjustwidth}
\end{table}
Table~\ref{tab:compare_time} shows that for all of the problem classes, the number of partial SVD required by \proc{ASALM} was more than twice the number that \proc{NSA} required. On the other hand, there was a big difference in CPU times; this difference can be explained by the fact that \proc{ASALM} required more leading singular values than \proc{NSA} did per partial SVD computation. Table~\ref{tab:compare_quality} shows that although the relative errors of the low-rank components produced by \proc{NSA} were slightly better, the relative errors of the sparse components produced by \proc{NSA} were significantly better than those produced by \proc{ASALM}. Finally, in Figure~\ref{fig:nsa_45dB}, we plot the decomposition of $D=X^0+S^0+\zeta^0\in\reals^{n\times n}$ generated by \proc{NSA}, where $\rank(X^0)=75$, $\norm{S^0}_0=112,500$ and ${\rm SNR}(D)=45$. In the first row, we plot randomly selected 1500 components of $S^0$ and 100 leading singular values of $X^0$ in the first row. In the second row, we plot the same components of $S^{sol}$ and 100 singular of $X^{sol}$ produced by \proc{NSA}. In the third row, we plot the absolute errors of $S^{sol}$ and $X^{sol}$. Note that the scales of the graphs showing absolute errors of $S^{sol}$ and $X^{sol}$ are larger than those of $S^0$ and $X^0$. And in the fourth row, we plot the same 1500 random components of $\zeta^0$. When we compare the absolute error graphs of $S^{sol}$ and $X^{sol}$ with the graph showing $\zeta^0$, we can confirm that the solution produced by \proc{NSA} is inline with Theorem~\ref{thm:candes2}.
\begin{table}[!htb]
    \begin{adjustwidth}{-2em}{-2em}
    \centering
    \caption{NSA vs ASALM: Solution time for decomposing $D\in\reals^{n\times n}$, $n=1500$}
    \renewcommand{\arraystretch}{1.3}
    {\footnotesize
    \begin{tabular}{ccc|c|c|c|c|}
    \cline{4-7}
    &&&$\mathbf{c_r}$=\textbf{0.05} $\mathbf{c_p}$=\textbf{0.05}&$\mathbf{c_r}$=\textbf{0.05} $\mathbf{c_p}$=\textbf{0.1}&$\mathbf{c_r}$=\textbf{0.1} $\mathbf{c_p}$=\textbf{0.05}&$\mathbf{c_r}$=\textbf{0.1} $\mathbf{c_p}$=\textbf{0.1}\\\hline
    \multicolumn{1}{|c|}{SNR}&\multicolumn{1}{|c|}{Alg.}& Field & min/\textbf{avg}/max & min/\textbf{avg}/max & min/\textbf{avg}/max & min/\textbf{avg}/max\\ \hline
    \multicolumn{1}{|c|}{\multirow{4}{*}{$80dB$}}
    &\multicolumn{1}{|c|}{\multirow{2}{*}{$\mathbf{NSA}$}}
    & $\mathbf{SVD\#}$ & 10/\textbf{10.0}/10 & 10/\textbf{10.9}/11 & 12/\textbf{12.0}/12 & 12/\textbf{12.2}/13 \\
    \multicolumn{1}{|c|}{}&\multicolumn{1}{|c|}{}& $\mathbf{CPU}$ & 38.6/\textbf{44.1}/46.6 & 43.7/\textbf{48.6}/51.9 & 78.6/\textbf{84.1}/90.8& 80.7/\textbf{97.7}/155.2 \\ \cline{2-7}
    \multicolumn{1}{|c|}{}
    &\multicolumn{1}{|c|}{\multirow{2}{*}{$\mathbf{ASALM}$}}
    & $\mathbf{SVD\#}$ & 22/\textbf{22.0}/22 & 20/\textbf{20.0}/20 & 29/\textbf{29.0}/29& 29/\textbf{29.4}/30 \\
    \multicolumn{1}{|c|}{}&\multicolumn{1}{|c|}{}& $\mathbf{CPU}$ & 657.3/\textbf{677.8}/736.2 & 809.7/\textbf{850.0}/874.7 & 1277.3/\textbf{1316.1}/1368.6	& 1833.2/\textbf{1905.2}/2004.7 \\ \hline
    \multicolumn{1}{|c|}{\multirow{4}{*}{$45dB$}}
    &\multicolumn{1}{|c|}{\multirow{2}{*}{$\mathbf{NSA}$}}
    & $\mathbf{SVD\#}$ & 7/\textbf{7.9}/8 & 8/\textbf{8.0}/8 & 9/\textbf{9.0}/9 & 9/\textbf{9.0}/9 \\
    \multicolumn{1}{|c|}{}&\multicolumn{1}{|c|}{}& $\mathbf{CPU}$ & 34.1/\textbf{37.7}/44.0 & 30.7/\textbf{37.1}/45.6 & 55.6/\textbf{59.0}/63.7& 55.9/\textbf{59.7}/64.8 \\ \cline{2-7}
    \multicolumn{1}{|c|}{}
    &\multicolumn{1}{|c|}{\multirow{2}{*}{$\mathbf{ASALM}$}}
    & $\mathbf{SVD\#}$ & 21/\textbf{21}/21 & 18/\textbf{18.5}/19 & 28/\textbf{28.0}/28& 27/\textbf{27.3}/28 \\
    \multicolumn{1}{|c|}{}&\multicolumn{1}{|c|}{}& $\mathbf{CPU}$ & 666.6/\textbf{686.9}/708.9 & 835.7/\textbf{857.1}/887.2 & 1201.9/\textbf{1223.2}/1277.5	& 1677.1/\textbf{1739.1}/1846.5 \\ \hline
    \end{tabular}
    \label{tab:compare_time}
    }
    \end{adjustwidth}
\end{table}
\begin{table}[!htb]
    \begin{adjustwidth}{-2em}{-2em}
    \centering
    \caption{NSA vs ASALM: Solution accuracy for decomposing $D\in\reals^{n\times n}$, $n=1500$}
    \renewcommand{\arraystretch}{1.5}
    {\footnotesize
    \begin{tabular}{ccc|c|c|c|c|}
    \cline{4-7}
    &&&$\mathbf{c_r}$=\textbf{0.05} $\mathbf{c_p}$=\textbf{0.05}&$\mathbf{c_r}$=\textbf{0.05} $\mathbf{c_p}$=\textbf{0.1}&$\mathbf{c_r}$=\textbf{0.1} $\mathbf{c_p}$=\textbf{0.05}&$\mathbf{c_r}$=\textbf{0.1} $\mathbf{c_p}$=\textbf{0.1}\\\hline
    \multicolumn{1}{|c|}{SNR}&\multicolumn{1}{|c|}{Alg.}&Relative Error& \textbf{avg}~/~max & \textbf{avg}~/~max & \textbf{avg}~/~max & \textbf{avg}~/~max\\ \hline
    \multicolumn{1}{|c|}{\multirow{4}{*}{$80dB$}}
    &\multicolumn{1}{|c|}{\multirow{2}{*}{$\mathbf{NSA}$}}
    & $\mathbf{\frac{\norm{X^{sol}-X^0}_F}{\norm{X^0}_F}}$
    &\textbf{1.8E-4}~/~2.2E-4&\textbf{2.1E-4}~/~2.6E-4&\textbf{1.3E-4}~/~1.3E-4&\textbf{2.8E-4}~/~2.9E-4\\
    \multicolumn{1}{|c|}{}&\multicolumn{1}{|c|}{}& $\mathbf{\frac{\norm{S^{sol}-S^0}_F}{\norm{S^0}_F}}$
    &\textbf{1.3E-4}~/~1.6E-4&\textbf{9.6E-5}~/~1.1E-4&\textbf{8.1E-5}~/~8.5E-5&\textbf{1.3E-4}~/~1.4E-4\\ \cline{2-7}
    \multicolumn{1}{|c|}{}
    &\multicolumn{1}{|c|}{\multirow{2}{*}{$\mathbf{ASALM}$}}
    & $\mathbf{\frac{\norm{X^{sol}-X^0}_F}{\norm{X^0}_F}}$
    &\textbf{3.9E-4}~/~4.2E-4&\textbf{8.4E-4}~/~8.8E-4&\textbf{6.6E-4}~/~6.8E-4&\textbf{1.4E-3}~/~1.4E-3\\
    \multicolumn{1}{|c|}{}&\multicolumn{1}{|c|}{}& $\mathbf{\frac{\norm{S^{sol}-S^0}_F}{\norm{S^0}_F}}$
    &\textbf{5.7E-4}~/~6.2E-4&\textbf{7.6E-4}~/~8.0E-4&\textbf{1.1E-3}~/~1.1E-3&\textbf{1.4E-3}~/~1.4E-3\\ \hline
    \multicolumn{1}{|c|}{\multirow{4}{*}{$45dB$}}
    &\multicolumn{1}{|c|}{\multirow{2}{*}{$\mathbf{NSA}$}}
    & $\mathbf{\frac{\norm{X^{sol}-X^0}_F}{\norm{X^0}_F}}$
    &\textbf{3.4E-3}~/~3.6E-3&\textbf{4.7E-3}~/~4.7E-3&\textbf{3.9E-3}~/~4.0E-3&\textbf{5.3E-3}~/~5.3E-3\\
    \multicolumn{1}{|c|}{}&\multicolumn{1}{|c|}{}& $\mathbf{\frac{\norm{S^{sol}-S^0}_F}{\norm{S^0}_F}}$
    &\textbf{1.8E-3}~/~1.8E-3&\textbf{2.3E-3}~/~2.3E-3&\textbf{2.6E-3}~/~3.5E-3&\textbf{3.1E-3}~/~3.1E-3\\ \cline{2-7}
    \multicolumn{1}{|c|}{}
    &\multicolumn{1}{|c|}{\multirow{2}{*}{$\mathbf{ASALM}$}}
    & $\mathbf{\frac{\norm{X^{sol}-X^0}_F}{\norm{X^0}_F}}$
    &\textbf{4.6E-3}~/~4.8E-3&\textbf{7.3E-3}~/~8.4E-3&\textbf{4.7E-3}~/~4.7E-3&\textbf{7.8E-3}~/~7.9E-3\\
    \multicolumn{1}{|c|}{}&\multicolumn{1}{|c|}{}& $\mathbf{\frac{\norm{S^{sol}-S^0}_F}{\norm{S^0}_F}}$
    &\textbf{4.8E-3}~/~4.9E-3&\textbf{5.8E-3}~/~7.0E-3&\textbf{5.5E-3}~/~5.5E-3&\textbf{7.3E-3}~/~7.5E-3\\ \hline
    \end{tabular}
    \label{tab:compare_quality}
    }
    \end{adjustwidth}
\end{table}
\subsection{Foreground Detection on a Noisy Video}
\label{sec:video_test_results}
We used \proc{NSA} and \proc{ASALM} to extract moving objects in an airport security video~\cite{Li04_1J}, which is a sequence of 201 grayscale frames of size $144 \times 176$. We assume that the airport security video~\cite{Li04_1J} was not corrupted by Gaussian noise. We formed the $i$-th column of the data matrix $D$ by stacking the columns of the $i^{th}$ frame into a long vector, i.e. $D$ is in $\reals^{25344\times 201}$. In order to have a noisy video with $\proc{SNR}=20dB$ signal-to-noise ratio~(\proc{SNR}), given $D$, we chose $\varrho = \norm{D}_F/(\sqrt{144\times 176\times 201}~10^{\proc{SNR}/20})$ and then obtained a noisy $D$ by $D = D+\varrho~randn(144*176,201)$, where $randn(m,n)$ produces a random matrix with independent standard Gaussian entries. Solving for $(X^*,S^*)=\argmin_{X,S\in\reals^{25344\times 201}}\{\norm{X}_*+\xi\norm{S}_1:\ \norm{X+S-D}_F\leq\delta\}$, we decompose $D$ into a low rank matrix $X^*$ and a sparse matrix $S^*$. We estimate the $i$-th frame background image with the $i$-th column of $X^*$ and estimate the $i$-th frame moving object with the $i$-th column of $S^*$. Both algorithms are terminated when $\frac{\norm{(X_{k+1},S_{k+1})-(X_{k},S_{k})}_F}{\norm{(X_{k},S_{k})}_F+1}\leq \varrho\times 10^{-4}$.

\begin{figure} [!h]
\centering
    \mbox{\hspace{4mm}$D(t)$:}
    \includegraphics[scale=1]{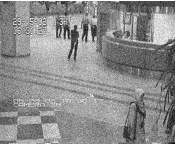}
    \includegraphics[scale=1]{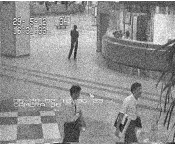}
    \includegraphics[scale=1]{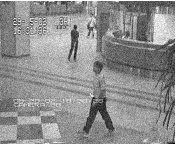}\\
    \mbox{$X^{sol}(t)$:}
    \includegraphics[scale=1]{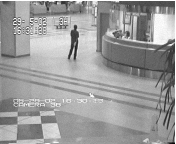}
    \includegraphics[scale=1]{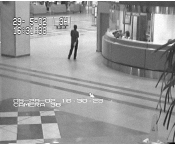}
    \includegraphics[scale=1]{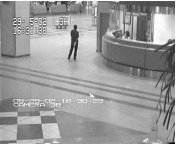}\\
    \mbox{$S^{sol}(t)$: }
    \includegraphics[scale=1]{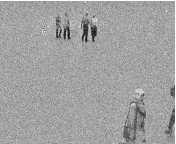}
    \includegraphics[scale=1]{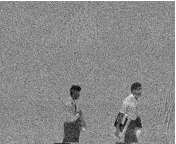}
    \includegraphics[scale=1]{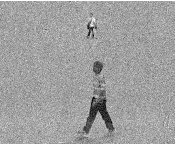}
\label{fig:noisy_reconstruction_test_nsa}
\caption{Background extraction from a video with 20dB SNR using \proc{NSA}}
\end{figure}

\begin{figure} [!h]
\centering
    \mbox{\hspace{4mm}$D(t)$:}
    \includegraphics[scale=1]{Dn_t35.png}
    \includegraphics[scale=1]{Dn_t100.png}
    \includegraphics[scale=1]{Dn_t125.png}\\
    \mbox{$X^{sol}(t)$:}
    \includegraphics[scale=1]{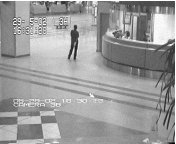}
    \includegraphics[scale=1]{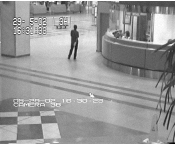}
    \includegraphics[scale=1]{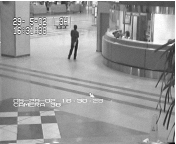}\\
    \mbox{$S^{sol}(t)$: }
    \includegraphics[scale=1]{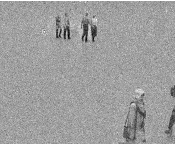}
    \includegraphics[scale=1]{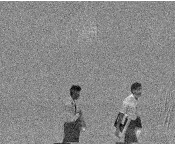}
    \includegraphics[scale=1]{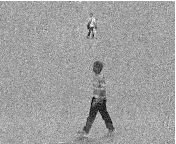}
\label{fig:noisy_reconstruction_test_asalm}
\caption{Background extraction from a video with 20dB SNR using \proc{ASALM}}
\end{figure}

The recovery statistics of each algorithm are are displayed in Table~\ref{tab:compare_video}. $(X^{sol},S^{sol})$ denote the variables corresponding to the low-rank and sparse components of $D$, respectively, when the algorithm of interest terminates. Figure~\ref{fig:noisy_reconstruction_test_nsa} and Figure~\ref{fig:noisy_reconstruction_test_asalm} show the $35$-th, $100$-th and $125$-th frames of the noise added airport security video~\cite{Li04_1J} in their first row of images. The second and third rows in these tables have the recovered background and foreground images of the selected frames, respectively. Even though the visual quality of recovered background and foreground are very similar, Table~\ref{tab:compare_video} shows that both the number of partial SVDs and the CPU time of \proc{NSA} are significantly less than those for \proc{ASALM}.
\begin{table}[!htb]
    \begin{adjustwidth}{-2em}{-2em}
    \centering
    \caption{NSA vs ASALM: Recovery statistics for foreground detection on a noisy video}
    \renewcommand{\arraystretch}{1.75}
    {\footnotesize
    \begin{tabular}{|c|c|c|c|c|c|c|}
    \hline
    Alg.& $\mathbf{CPU}$ & $\mathbf{SVD\#}$ & $\mathbf{\norm{X^{sol}}_*}$ & $\mathbf{\norm{S^{sol}}_1}$ & $\Rank(X^{sol})$ & $\mathbf{\frac{\norm{X^{sol}+S^{sol}-D}_F}{\norm{D}_F}}$ \\ \hline
    $\mathbf{NSA}$ & 160.8 & 19 & 398662.9 & 76221854.1 & 81 & 0.00068 \\ \hline
    $\mathbf{ASALM}$ & 910.0 & 94 & 401863.6 & 75751977.1 & 89 & 0.00080 \\ \hline
    \end{tabular}
    \label{tab:compare_video}
    }
    \end{adjustwidth}
\end{table}
\section{Acknowledgements}
We would like to thank to Min Tao for providing the code \proc{ASALM}.
\clearpage
\bibliographystyle{siam}
\bibliography{thesis}
\appendix
\section{Proof of Theorem~\ref{thm:alm}}
\label{sec:thm_proof}
\begin{definition}
Let $\phi:\reals^{m\times n}\rightarrow\reals$ and $\psi:\reals^{m\times n}\times\reals^{m\times n}\rightarrow\reals$ be closed convex functions and define
\begin{align}
Q^\phi(Z,S|X):=\psi(Z,S)+\phi(X)+\fprod{\gamma^\phi(X), Z-X}+\frac{\rho}{2}\norm{Z-X}_F^2,\\
Q^\psi(Z|X,S):=\phi(Z)+\psi(X,S)+\fprod{\gamma_x^\psi(X,S), Z-X}+\frac{\rho}{2}\norm{Z-X}_F^2, \label{eq:Q_psi}
\end{align}
and
\begin{align}
\left(p^\phi_x(X),p^\phi_s(X)\right):=\argmin_{Z,S\in\reals^{m\times n}}Q^\phi(Z,S|X), \label{eq:optimal_point_f}\\
p^\psi(X,S):=\argmin_{Z\in\reals^{m\times n}}Q^\psi(Z|X,S), \label{eq:optimal_point_psi}
\end{align}
where $\gamma^\phi(X)$ is any subgradient in the subdifferential $\partial \phi$ at the point $X$ and $\Big(\gamma_x^\psi(X,S),\gamma_s^\psi(X,S)\Big)$ is any subgradient in the subdifferential $\partial \psi$ at the point $(X,S)$.
\end{definition}
\begin{lemma}
\label{lem:alms}
Let $\phi$, $\psi$, $Q^\phi$, $Q^\psi$, $p_x^\phi$, $p_s^\phi$, $p^\psi$, $\gamma^\phi$, $\gamma_x^\psi$, $\gamma_s^\psi$ be as given in Definition A.1. and $\Phi(X,S):=\phi(X)+\psi(X,S)$. Let $X^0\in\reals^{m\times n}$ and define $\hat{X} := p_x^\phi(X^0)$ and $\hat{S} := p_s^\phi(X^0)$. If
\begin{align}
\Phi(\hat{X},\hat{S})\leq Q^\phi(\hat{X},\hat{S}|X^0), \label{eq:upper_envelope_condition_1}
\end{align}
then for any $(X,S)\in\reals^{m\times n}\times\reals^{m\times n}$,
\begin{align}
\frac{2}{\rho}\left(\Phi(X,S)-\Phi\left(\hat{X},\hat{S}\right)\right)\geq \norm{X-\hat{X}}_F^2-\norm{X-X^0}_F^2.\label{eq:f_lemma}
\end{align}
Moreover, if
\begin{align}
\Phi\left(p^\psi\left(\hat{X},\hat{S}\right),\hat{S}\right)\leq Q^\psi\left(p^\psi\left(\hat{X},\hat{S}\right)\Big|~\hat{X},\hat{S}\right), \label{eq:upper_envelope_condition_2}
\end{align}
then for any $(X,S)\in\reals^{m\times n}\times\reals^{m\times n}$,
\begin{align}
\frac{2}{\rho}\Bigg(\Phi(X,S)-\Phi\Big(p^\psi\left(\hat{X},\hat{S}\right),\hat{S}\Big)\Bigg)\geq \norm{X-p^\psi\left(\hat{X},\hat{S}\right)}_F^2-\norm{X-\hat{X}}_F^2.\label{eq:psi_lemma}
\end{align}
\end{lemma}
\begin{proof}
Let $X^0\in\reals^{m\times n}$ satisfy \eqref{eq:upper_envelope_condition_1}. Then for any $(X,S)\in\reals^{m\times n}\times\reals^{m\times n}$, we have
\begin{align}
\Phi(X,S)-\Phi\left(p^\phi_x(X^0),\hat{S}\right)&\geq \Phi(X,S)-Q^\phi\left(\hat{X},\hat{S}|X^0\right).
\end{align}
First order optimality conditions for \eqref{eq:optimal_point_f} and $\psi$ being a closed convex function guarantee that there exists $\Bigg(\gamma^\psi_x\Big(\hat{X},\hat{S}\Big),\gamma^\psi_s\Big(\hat{X},\hat{S}\Big)\Bigg)\in\partial\psi\Big(\hat{X},\hat{S}\Big)$ such that
\begin{align}
\gamma^\psi_x\Big(\hat{X},\hat{S}\Big)+\gamma^\phi(X^0)+\rho\Big(\hat{X}-X^0\Big)=0,\label{eq:psi_x_subgradient_1}\\
\gamma^\psi_s\Big(\hat{X},\hat{S}\Big)=0,\label{eq:psi_s_subgradient_1}
\end{align}
where $\partial\psi\Big(\hat{X},\hat{S}\Big)$ denotes the subdifferential of $\psi(.,.)$ at the point $\Big(\hat{X},\hat{S}\Big)$.

Moreover, using the convexity of $\psi(.,.)$ and $\phi(.)$, we have
\begin{align*}
&\psi(X,S)\geq\psi\Big(\hat{X},\hat{S}\Big)+\Big\langle\gamma^\psi_x\Big(\hat{X},\hat{S}\Big), X-\hat{X}\Big\rangle+\Big\langle\gamma^\psi_s\Big(\hat{X},\hat{S}\Big), S-\hat{S}\Big\rangle,\\
&\phi(X)\geq \phi(X^0)+\fprod{\gamma^\phi(X^0), X-X^0}.
\end{align*}
These two inequalities and \eqref{eq:psi_s_subgradient_1} together imply
\begin{align}
\Phi(X,S)\geq \psi\Big(\hat{X},\hat{S}\Big)+\Big\langle\gamma^\psi_x\Big(\hat{X},\hat{S}\Big), X-\hat{X}\Big\rangle+\phi(X^0)+\Big\langle\gamma^\phi(X^0), X-X^0\Big\rangle.
\end{align}
This inequality together with \eqref{eq:upper_envelope_condition_1} and \eqref{eq:psi_x_subgradient_1} gives
\begin{align*}
&\Phi(X,S)-\Phi\Big(\hat{X},\hat{S}\Big)\\
\geq &\Big\langle\gamma^\psi_x\Big(\hat{X},\hat{S}\Big), X-\hat{X}\Big\rangle+\Big\langle\gamma^\phi(X^0), X-X^0\Big\rangle-\Big\langle\gamma^\phi(X^0), \hat{X}-X^0\Big\rangle-\frac{\rho}{2}\norm{X-X^0}_F^2,\\
= &\Big\langle\gamma^\phi(X^0)+\gamma^\psi_x\Big(\hat{X},\hat{S}\Big),\ X-\hat{X}\Big\rangle-\frac{\rho}{2}\norm{X-X^0}_F^2,\\
= &\rho~\Big\langle X^0-\hat{X},\ X-\hat{X}\Big\rangle-\frac{\rho}{2}\norm{X-X^0}_F^2,\\
= &\frac{\rho}{2}\Big(\norm{X-\hat{X}}_F^2-\norm{X-X^0}_F^2\Big).
\end{align*}
Hence, we have \eqref{eq:f_lemma}. Suppose that $X^0$ satisfies \eqref{eq:upper_envelope_condition_2}. Then for any $(X,S)\in\reals^{m\times n}\times\reals^{m\times n}$, we have
\begin{align}
\Phi(X,S)-\Phi\Big(p^\psi\left(\hat{X},\hat{S}\right),\hat{S}\Big)&\geq \Phi(X,S)-Q^\psi\left(p^\psi\left(\hat{X},\hat{S}\right)|~\hat{X},\hat{S}\right).
\end{align}
First order optimality conditions for \eqref{eq:optimal_point_psi} and $\phi$ being a closed convex function guarantee that there exists $\gamma^\phi\Big(p^\psi\left(\hat{X},\hat{S}\right)\Big)\in\partial \phi\Big(p^\psi\left(\hat{X},\hat{S}\right)\Big)$ such that
\begin{align}
\gamma^\phi\Big(p^\psi\left(\hat{X},\hat{S}\right)\Big)+\gamma_x^\psi\Big(\hat{X},\hat{S}\Big)+\rho\Big(p^\psi\left(\hat{X},\hat{S}\right)- \hat{X}\Big)=0.\label{eq:f_subgradient_1}
\end{align}
Moreover, using the convexity of $\phi(.)$ and $\psi(.,.)$, we have
\begin{align}
\phi(X)\geq \phi\Big(p^\psi\left(\hat{X},\hat{S}\right)\Big)+\Big\langle\gamma^\phi\Big(p^\psi\left(\hat{X},\hat{S}\right)\Big),X-p^\psi\left(\hat{X},\hat{S}\right)\Big\rangle, \label{eq:f_convexity_2}\\
\psi(X,S)\geq\psi\Big(\hat{X},\hat{S}\Big)+\Big\langle\gamma^\psi_x\left(\hat{X},\hat{S}\right), X-\hat{X}\Big\rangle, \label{eq:psi_convexity_2}
\end{align}
where \eqref{eq:psi_convexity_2} follows from the fact that $\left(\hat{X},\hat{S}\right)=\argmin_{X,S}Q^\phi(X,S|~X^0)$ implies $\left(\gamma^\psi_x\left(\hat{X},\hat{S}\right),~0\right)\in\partial \psi\left(\hat{X},\hat{S}\right)$, i.e. we can set $\gamma^\psi_s\left(\hat{X},\hat{S}\right)=0$.
Summing the two inequalities \eqref{eq:f_convexity_2} and \eqref{eq:psi_convexity_2} give
\begin{align}
\Phi(X,S)\geq &\ \psi\left(\hat{X},\hat{S}\right)+\Big\langle\gamma^\psi_x\left(\hat{X},\hat{S}\right), X-\hat{X}\Big\rangle+\phi\left(p^\psi\left(\hat{X},\hat{S}\right)\right)+\Big\langle\gamma^\phi\left(p^\psi\left(\hat{X},\hat{S}\right)\right),X-p^\psi\left(\hat{X},\hat{S}\right)\Big\rangle.
\end{align}
This inequality together with \eqref{eq:upper_envelope_condition_2} and \eqref{eq:f_subgradient_1} gives
\begin{align*}
&\Phi(X,S)-\Phi\Big(p^\psi\left(\hat{X},\hat{S}\right),\hat{S}\Big)\\
\geq &\Big\langle\gamma^\psi_x(\hat{X},\hat{S}), X-\hat{X}\Big\rangle+\Big\langle\gamma^\phi\Big(p^\psi\left(\hat{X},\hat{S}\right)\Big),~X-p^\psi\left(\hat{X},\hat{S}\right)\Big\rangle\\
&-\Big\langle\gamma_x^\psi\Big(\hat{X},\hat{S}\Big),~p^\psi\Big(\hat{X},\hat{S}\Big)-\hat{X}\Big\rangle
-\frac{\rho}{2}\norm{p^\psi\Big(\hat{X},\hat{S}\Big)-\hat{X}}_F^2,\\
= &\Big\langle\gamma^\phi\Big(p^\psi\left(\hat{X},\hat{S}\right)\Big)+\gamma^\psi_x\Big(\hat{X},\hat{S}\Big),\ X-p^\psi\left(\hat{X},\hat{S}\right)\Big\rangle-\frac{\rho}{2}\norm{p^\psi\left(\hat{X},\hat{S}\right)-\hat{X}}_F^2,\\
= &\rho~\Big\langle \hat{X}-p^\psi\left(\hat{X},\hat{S}\right),\ X-p^\psi\left(\hat{X},\hat{S}\right)\Big\rangle-\frac{\rho}{2}\norm{p^\psi\left(\hat{X},\hat{S}\right)-\hat{X}}_F^2,\\
= &\frac{\rho}{2}\left(\norm{X-p^\psi\left(\hat{X},\hat{S}\right)}_F^2-\norm{X-\hat{X}}_F^2\right).
\end{align*}
Hence, we have \eqref{eq:psi_lemma}.
\end{proof}

We are now ready to give the proof of Theorem~\ref{thm:alm}.
\begin{proof}
Let $I:=\{0\leq i\leq k-1: \Phi(X_{i+1},S_i)\leq\mathcal{L}_\rho(X_{i+1},Z_i,S_i;Y_i)\}$ and $I^c:=\{0,1,...,k-1\}\setminus I$. Since $\grad \phi(.)$ is Lipschitz continuous with Lipschitz constant $L$ and $\rho\geq L$, $\Phi(p^\phi_x(X),p^\phi_s(X))\leq Q^\phi(p^\phi_x(X),p^\phi_s(X)|~X)$ is true for all $X\in\reals^{m\times n}$. Since \eqref{eq:upper_envelope_condition_1} in Lemma~\ref{lem:alms} is true for all $X^0\in\reals^{m\times n}$, \eqref{eq:f_lemma} is true for all $(X,S)\in\reals^{m\times n}\times\reals^{m\times n}$. Particularly, since for all $i\in I\cup I^c$
\begin{align}
(Z_{i+1},S_{i+1})=\argmin_{Z,S} Q^\phi(Z,S|~X_{i+1}),\label{eq:zs_optimality}
\end{align}
setting $(X,S):=(X^*,S^*)$ and $X^0:=X_{i+1}$ in Lemma~\ref{lem:alms} imply that $p^\phi_x(X_{i+1})=Z_{i+1}$, $p^\phi_s(X_{i+1})=S_{i+1}$ and we have
\begin{align}
\frac{2}{\rho}\left(\Phi(X^*,S^*)-\Phi(Z_{i+1},S_{i+1})\right)\geq \norm{Z_{i+1}-X^*}_F^2-\norm{X_{i+1}-X^*}_F^2.\label{eq:ineq1}
\end{align}
Moreover, \eqref{eq:zs_optimality} implies that for all $i\in I\cup I^c$, there exits $\left(\gamma^\psi_x(Z_i,S_i),\gamma^\psi_s(Z_i,S_i)\right)\in\partial\psi(Z_i,S_i)$ such that
\begin{align}
\gamma^\psi_x(Z_i,S_i)+\grad \phi(X_i)+ \rho(Z_i-X_i)=0, \label{eq:subgrad_psi_x}\\
\gamma^\psi_s(Z_i,S_i)=0.
\end{align}
\eqref{eq:subgrad_psi_x} and the definition of $Y_{i+1}$ of Algorithm~\proc{ALM-S} shown in \textbf{Algorithm~\ref{alg:alms}} imply that
$$\gamma^\psi_x(Z_i,S_i)=-\grad \phi(X_i)+\rho(X_i-Z_i)=Y_i.$$
Hence, by defining $Q^\psi(.|~Z_i,S_i)$ according to \eqref{eq:Q_psi} using $\gamma^\psi_x(Z_i,S_i)=Y_i$, for all $X\in\reals^{m\times n}$ we have
\begin{align}
&\mathcal{L}_\rho(X,Z_i,S_i;Y_i)= \phi(X)+\psi(Z_i,S_i)+\fprod{Y_i, X-Z_i}+\frac{\rho}{2}\norm{X-Z_i}_F^2= Q^\psi(X|~Z_i,S_i).
\end{align}
for all $i\in I\cup I^c$.
Hence, for all $i\in I$ $X_{i+1} =\argmin_X \mathcal{L}_\rho(X,Z_i,S_i;Y_i)=\argmin_X Q^\psi(X|Z_i,S_i)$. Thus, for all $i\in I$, setting $X^0:=X_i$ in Lemma~\ref{lem:alms} imply $p^\phi_x(X_i)=Z_i$, $p^\phi_s(X_i)=S_i$ and $p^\psi(p^\phi_x(X_i),p^\phi_s(X_i))=p^\psi(Z_i,S_i)=X_{i+1}$. For all $i\in I$ we have $\Phi(X_{i+1},S_i)\leq\mathcal{L}_\rho(X_{i+1},Z_i,S_i;Y_i)=Q^\psi(X_{i+1}|Z_i,S_i)$. Hence, for all $i\in I$ setting $X^0:=X_i$ in Lemma~\ref{lem:alms} satisfies \eqref{eq:upper_envelope_condition_2}. Therefore, setting $(X,S):=(X^*,S^*)$ and $X^0:=X_{i}$ in Lemma~\ref{lem:alms} implies that
\begin{align}
\frac{2}{\rho}\left(\Phi(X^*,S^*)-\Phi(X_{i+1},S_i)\right)\geq \norm{X_{i+1}-X^*}_F^2-\norm{Z_i-X^*}_F^2.\label{eq:ineq2}
\end{align}
For any $i\in I$, summing \eqref{eq:ineq1} and \eqref{eq:ineq2} gives
\begin{align}
\frac{2}{\rho}\left(2\Phi(X^*,S^*)-\Phi(X_{i+1},S_i)-\Phi(Z_{i+1},S_{i+1})\right)\geq \norm{Z_{i+1}-X^*}_F^2-\norm{Z_i-X^*}_F^2.\label{eq:I_ineq}
\end{align}
Moreover, since $X_{i+1}=Z_i$ for $i\in I^c$ and \eqref{eq:ineq1} holds for all $i\in I\cup I^c$, we trivially have
\begin{align}
\frac{2}{\rho}\left(\Phi(X^*,S^*)-\Phi(Z_{i+1},S_{i+1})\right)\geq \norm{Z_{i+1}-X^*}_F^2-\norm{Z_i-X^*}_F^2.\label{eq:Ic_ineq}
\end{align}
Summing \eqref{eq:I_ineq} and \eqref{eq:Ic_ineq} over $i=0,1,...,k-1$ gives
\begin{align}
\frac{2}{\rho}\left(\left(2|I|+|I^c|\right)\Phi(X^*,S^*)-\sum_{i\in I}\Phi(X_{i+1},S_i)-\sum_{i=0}^{k-1}\Phi(Z_{i+1},S_{i+1})\right)\geq \norm{Z_{k}-X^*}_F^2-\norm{Z_0-X^*}_F^2.\label{eq:key_ineq}
\end{align}
For any $i\in I\cup I^c$, setting $(X,S):=(X_{i+1},S_{i})$ and $X^0:=X_{i+1}$ in Lemma~\ref{lem:alms} gives
\begin{align}
\frac{2}{\rho}\left(\Phi(X_{i+1},S_{i})-\Phi(Z_{i+1},S_{i+1})\right)\geq \norm{Z_{i+1}-X_{i+1}}_F^2\geq 0.\label{eq:ineq3}
\end{align}
Trivially, for $i=1,...,k$ we also have
\begin{align}
\frac{2}{\rho}\left(\Phi(X_{i},S_{i-1})-\Phi(Z_{i},S_{i})\right)\geq \norm{Z_{i}-X_{i}}_F^2\geq 0.\label{eq:ineq4}
\end{align}
Moreover, since for all $i\in I$ setting $X^0:=X_i$ in Lemma~\ref{lem:alms} satisfies \eqref{eq:upper_envelope_condition_2}, setting $(X,S):=(Z_i,S_i)$ and $X^0:=X_{i}$ in Lemma~\ref{lem:alms} implies that
\begin{align}
\frac{2}{\rho}\left(\Phi(Z_i,S_i)-\Phi(X_{i+1},S_i)\right)\geq \norm{X_{i+1}-Z_i}_F^2\geq 0.\label{eq:ineq5_I}
\end{align}
And since $X_{i+1}=Z_i$ for all $i\in I^c$, \eqref{eq:ineq5_I} trivially holds for all $i\in I^c$. Thus, for all $i\in I\cup I^c$ we have
\begin{align}
\frac{2}{\rho}\left(\Phi(Z_i,S_i)-\Phi(X_{i+1},S_i)\right)\geq 0.\label{eq:ineq5}
\end{align}
Adding \eqref{eq:ineq3} and \eqref{eq:ineq5} yields $\Phi(Z_i,S_i)\geq \Phi(Z_{i+1},S_{i+1})$ for all $i\in I\cup I^c$ and adding \eqref{eq:ineq4} and \eqref{eq:ineq5} yields $\Phi(X_i,S_{i-1})\geq \Phi(X_{i+1},S_i)$ for all $i=1,...,k-1$. Hence,
\begin{align}
\sum_{i=0}^{k-1}\Phi(Z_{i+1},S_{i+1})\geq k \Phi(Z_{k},S_{k}), \mbox{ and } \sum_{i\in I}\Phi(X_{i+1},S_i) \geq n_k \Phi(X_{k},S_{k-1}).
\end{align}
These two inequalities, \eqref{eq:key_ineq} and the fact that $X_0=Z_0$ imply
\begin{align}
\frac{2}{\rho}\left(\left(2|I|+|I^c|\right)\Phi(X^*,S^*)-n_k \Phi(X_{k},S_{k-1})-k \Phi(Z_{k},S_{k})\right)\geq -\norm{X_0-X^*}_F^2.
\end{align}
Hence, \eqref{eq:alms_thm} follows from the facts: $2|I|+|I^c|=k+n_k$ and $n_k \Phi(X_{k},S_{k-1})+k \Phi(Z_{k},S_{k})\geq (k+n_k)\Phi(Z_{k},S_{k})$ due to \eqref{eq:ineq3}.
\end{proof}
\section{Proof of Lemma~\ref{lem:finite_sums}}
\label{sec:lem_proof}
\begin{proof}
Since $Y^*$ and $\theta^*$ are optimal Lagrangian dual variables, we have
\begin{align*}
(X^*,X^*,S^*)=\argmin_{X,Z,S}\norm{X}_*+\xi~\norm{S}_1+\fprod{Y^*, X-Z}+\frac{\theta^*}{2}\left( \norm{Z+S-D}^2_F-\delta^2\right).
\end{align*}
Then from first-order optimality conditions, we have
\begin{align*}
&0\in\partial\norm{X^*}_*+Y^*,\\
&0\in\xi~\partial\norm{S^*}_1+\theta^*(X^*+S^*-D),\\
&-Y^*+\theta^*(X^*+S^*-D)=0.
\end{align*}
Hence, $-Y^*\in\partial\norm{X^*}_*$ and $-Y^*\in\xi~\partial\norm{S^*}_1$.

For $k\geq 0$,
since $X_{k+1}$ is the optimal solution for the $k$-th subproblem given in Step~\ref{algeq:subproblem1} in \textbf{Algorithm~\ref{alg:nsa}}, from the first-order optimality conditions it follows that
\begin{align}
\label{eq:Xopt_cond}
0\in\partial\norm{X_{k+1}}_*+ Y_k+\rho_k(X_{k+1}-Z_k).
\end{align}
For $k\geq 0$,
let $\theta_k\geq 0$ be the optimal Lagrange multiplier for the quadratic constraint in the $k$-th subproblem given in Step~\ref{algeq:subproblem2} in \textbf{Algorithm~\ref{alg:nsa}}. Since $(S_{k+1},Z_{k+1})$ is the optimal solution, from the first-order optimality conditions it follows that
\begin{align}
0\in\xi\partial\norm{S_{k+1}}_1+ \theta_k(Z_{k+1}+S_{k+1}-D), \label{eq:Sopt_cond}\\
-Y_k+\rho_k(Z_{k+1}-X_{k+1})+\theta_k(Z_{k+1}+S_{k+1}-D)=0. \label{eq:Zopt_cond}
\end{align}
From \eqref{eq:Xopt_cond}, it follows that $-\hat{Y}_{k+1}\in\partial\norm{X_{k+1}}_*$. Hence, $\{\hat{Y}_k\}_{k\in\integers_+}$ is a bounded sequence. From \eqref{eq:Sopt_cond} and \eqref{eq:Zopt_cond}, it follows that $-Y_{k+1}\in\xi~\partial\norm{S_{k+1}}_1$. Hence, $\{Y_k\}_{k\in\integers_+}$ is also a bounded sequence.

Furthermore, since $Y_{k+1}-Y_k=\rho_k(X_{k+1}-Z_{k+1})$ and $Y_{k+1}-\hat{Y}_{k+1}=\rho_k(Z_k-Z_{k+1})$, we have
\begin{align*}
 &\rho_k^{-1}\fprod{Y_{k+1}-Y_k, Y_{k+1}-Y^*}\\
=&\fprod{X_{k+1}-Z_{k+1}, Y_{k+1}-Y^*},\\
=&\fprod{X_{k+1}-X^*, Y_{k+1}-Y^*}+\fprod{X^*-Z_{k+1}, Y_{k+1}-Y^*},\\
=&\fprod{X_{k+1}-X^*, Y_{k+1}-\hat{Y}_{k+1}}+\fprod{X_{k+1}-X^*, \hat{Y}_{k+1}-Y^*}+\fprod{X^*-Z_{k+1}, Y_{k+1}-Y^*},\\
=&\rho_k\fprod{X_{k+1}-X^*, Z_k-Z_{k+1}}+\fprod{X_{k+1}-X^*, \hat{Y}_{k+1}-Y^*}+\fprod{X^*-Z_{k+1}, Y_{k+1}-Y^*}.
\end{align*}
Using the above equality, for all $k\geq 1$, we trivially have
\begin{align}
&\norm{Z_{k+1}-X^*}_F^2+\rho_{k}^{-2}\norm{Y_{k+1}-Y^*}_F^2 \nonumber\\
=&\norm{Z_{k}-X^*}_F^2+\rho_{k}^{-2}\norm{Y_{k}-Y^*}_F^2-\norm{Z_{k+1}-Z_k}_F^2-\rho_{k}^{-2}\norm{Y_{k+1}-Y_k}_F^2, \nonumber\\
&+2\fprod{Z_{k+1}-X^*, Z_{k+1}-Z_k}+2\rho_k^{-2}\fprod{Y_{k+1}-Y_k, Y_{k+1}-Y^*}, \nonumber\\
=&\norm{Z_{k}-X^*}_F^2+\rho_{k}^{-2}\norm{Y_{k}-Y^*}_F^2-\norm{Z_{k+1}-Z_k}_F^2-\rho_{k}^{-2}\norm{Y_{k+1}-Y_k}_F^2, \nonumber\\
&+2\fprod{Z_{k+1}-X^*, Z_{k+1}-Z_k}+2\fprod{X_{k+1}-X^*, Z_k-Z_{k+1}} \nonumber\\
&-2\rho_k^{-1}\left(\fprod{-\hat{Y}_{k+1}+Y^*, X_{k+1}-X^*}+\fprod{-Y_{k+1}+Y^*, X^*-Z_{k+1}}\right) \nonumber\\
=&\norm{Z_{k}-X^*}_F^2+\rho_{k}^{-2}\norm{Y_{k}-Y^*}_F^2-\norm{Z_{k+1}-Z_k}_F^2-\rho_{k}^{-2}\norm{Y_{k+1}-Y_k}_F^2,\nonumber\\
&+2\fprod{Z_{k+1}-X_{k+1}, Z_{k+1}-Z_k}-2\rho_k^{-1}\left(\fprod{-\hat{Y}_{k+1}+Y^*, X_{k+1}-X^*}+\fprod{-Y_{k+1}+Y^*, X^*-Z_{k+1}}\right)\nonumber\\
=&\norm{Z_{k}-X^*}_F^2+\rho_{k}^{-2}\norm{Y_{k}-Y^*}_F^2-\norm{Z_{k+1}-Z_k}_F^2-\rho_{k}^{-2}\norm{Y_{k+1}-Y_k}_F^2, \nonumber\\
&-2\rho_k^{-1}\left(\fprod{Y_{k+1}-Y_k, Z_{k+1}-Z_k}+\fprod{-\hat{Y}_{k+1}+Y^*, X_{k+1}-X^*}+\fprod{-Y_{k+1}+Y^*, X^*-Z_{k+1}}\right) \label{eq:preinduction_step}
\end{align}
Since $-Y_k\in\xi~\partial\norm{S_k}_1$ for all $k\geq 1$ and $-Y^*\in\xi~\partial\norm{S^*}_1$, we have for all $k\geq 1$
\begin{align}
\fprod{-Y_{k+1}+Y_k, S_{k+1}-S_k}\geq 0, \label{eq:S_subgradient_monotone}\\
\fprod{-Y_{k+1}+Y^*, S_{k+1}-S^*}\geq 0. \label{eq:Sopt_subgradient_monotone}
\end{align}
Since $\rho_{k+1}\geq\rho_k$ for all $k\geq 1$, adding \eqref{eq:S_subgradient_monotone}, \eqref{eq:Sopt_subgradient_monotone} to \eqref{eq:preinduction_step} and subtracting \eqref{eq:Sopt_subgradient_monotone} from \eqref{eq:preinduction_step}, we have
\begin{align}
&\norm{Z_{k+1}-X^*}_F^2+\rho_{k+1}^{-2}\norm{Y_{k+1}-Y^*}_F^2 \nonumber\\
\leq~&\norm{Z_{k+1}-X^*}_F^2+\rho_{k}^{-2}\norm{Y_{k+1}-Y^*}_F^2 \nonumber\\
\leq~&\norm{Z_{k}-X^*}_F^2+\rho_{k}^{-2}\norm{Y_{k}-Y^*}_F^2-\norm{Z_{k+1}-Z_k}_F^2-\rho_{k}^{-2}\norm{Y_{k+1}-Y_k}_F^2 \nonumber\\
&-2\rho_k^{-1}\left(\fprod{-\hat{Y}_{k+1}+Y^*, X_{k+1}-X^*}+\fprod{-Y_{k+1}+Y^*, S_{k+1}-S^*}\right) \nonumber\\
&-2\rho_k^{-1}\left(\fprod{Y_{k+1}-Y_k, Z_{k+1}+S_{k+1}-Z_k-S_k}+\fprod{-Y_{k+1}+Y^*, X^*+S^*-Z_{k+1}-S_{k+1}}\right) \label{eq:induction_step}
\end{align}
Applying Lemma~\ref{lem:chi_subgradient} on the $k$-th subproblem given in Step~\ref{algeq:subproblem2} in Algorithm~\ref{alg:nsa}, it follows that
$$(Y_{k+1},Y_{k+1})\in\partial \mathbf{1}_{\chi}(Z_{k+1},S_{k+1}).$$
Using arguments similar to those used in the proof of Lemma~\ref{lem:chi_subgradient}, one can also show that
$$(Y^*,Y^*)\in\partial \mathbf{1}_{\chi}(X^*,S^*).$$
Moreover, since $-Y_{k}\in\xi~\partial\norm{S_k}_1$, $-\hat{Y}_{k}\in\xi~\partial\norm{X_k}_*$ for all $k\geq 1$, $-Y^*\in\xi~\partial\norm{S^*}_1$ and $-Y^*\in\partial\norm{X^*}_*$, we have
\begin{align*}
\fprod{Y_{k+1}-Y_k, Z_{k+1}+S_{k+1}-Z_k-S_k}\geq 0,\\
\fprod{-Y_{k+1}+Y^*, X^*+S^*-Z_{k+1}-S_{k+1}}\geq 0,\\
\fprod{-Y_{k+1}+Y^*, S_{k+1}-S^*}\geq 0,\\
\fprod{-\hat{Y}_{k+1}+Y^*, X_{k+1}-X^*}\geq 0,
\end{align*}
for all $k\geq 1$.
Therefore, the above inequalities and \eqref{eq:induction_step} together imply that $\{\norm{Z_{k}-X^*}_F^2+\rho_{k}^{-2}\norm{Y_{k}-Y^*}_F^2\}_{k\in\integers_+}$ is a non-increasing sequence. Moreover, we also have
\begin{align*}
&\sum_{k\in\integers_+}\norm{Z_{k+1}-Z_k}_F^2+\rho_{k}^{-2}\norm{Y_{k+1}-Y_k}_F^2 \\
+&2\sum_{k\in\integers_+}\rho_k^{-1}\left(\fprod{-\hat{Y}_{k+1}+Y^*, X_{k+1}-X^*}+\fprod{-Y_{k+1}+Y^*, S_{k+1}-S^*}\right)\\
+&2\sum_{k\in\integers_+}\rho_k^{-1}\left(\fprod{Y_{k+1}-Y_k, Z_{k+1}+S_{k+1}-Z_k-S_k}+\fprod{-Y_{k+1}+Y^*, X^*+S^*-Z_{k+1}-S_{k+1}}\right)\\
=&\sum_{k\in\integers_+}\left(\norm{Z_{k}-X^*}_F^2+\rho_{k}^{-2}\norm{Y_{k}-Y^*}_F^2-\norm{Z_{k+1}-X^*}_F^2-\rho_{k+1}^{-2}\norm{Y_{k+1}-Y^*}_F^2\right)<\infty
\end{align*}
\end{proof}
\newpage
\section{Additional Statistics for Numerical Experiments}
\begin{table}[!htb]
    \begin{adjustwidth}{-2em}{-2em}
    \centering
    \caption{NSA vs ASALM: Additional statistics on solution accuracy for decomposing $D\in\reals^{n\times n}$, $n=1500$, ${\rm SNR}(D)=80dB$}
    \renewcommand{\arraystretch}{1.4}
    {\footnotesize
    \begin{tabular}{cc|c|c|c|}
    \cline{2-5}
    &\multicolumn{2}{|c|}{$\mathbf{c_r}$=\textbf{0.05} $\mathbf{c_p}$=\textbf{0.05}}&\multicolumn{2}{|c|}{$\mathbf{c_r}$=\textbf{0.05} $\mathbf{c_p}$=\textbf{0.1}}\\ \cline{2-5}
    &\multicolumn{1}{|c|}{\textbf{NSA}} & \textbf{ASALM} & \textbf{NSA} & \textbf{ASALM} \\ \hline
    \multicolumn{1}{|c|}{\textbf{Error Type}}&\multicolumn{1}{|c|}{\textbf{avg}~/~max}& \textbf{avg}~/~max & \textbf{avg}~/~max & \textbf{avg}~/~max\\ \hline
    \multicolumn{1}{|c|}{$\mathbf{|~\norm{X^{sol}}_*-\norm{X^0}_*|/\norm{X^0}_*}$}
    &\textbf{1.7E-6}~/~    5.2E-6&    \textbf{6.9E-6}~/~ 1.0E-5&    \textbf{5.0E-6}~/~ 3.2E-5&    \textbf{2.3E-5}~/~ 3.8E-5\\ \hline
    \multicolumn{1}{|c|}{$\mathbf{\max\{|\sigma_i-\sigma^0_i|: \sigma^0_i>0\}}$}
    &\textbf{4.1E-2}~/~	5.2E-2&	\textbf{3.9E-2}~/~	4.7E-2&	\textbf{5.2E-2}~/~ 1.8E-1&	\textbf{1.1E-1}~/~	1.7E-1\\ \hline
    \multicolumn{1}{|c|}{$\mathbf{\max\{|\sigma_i|: \sigma^0_i=0\}}$}
    &\textbf{7.9E-13}~/~	2.2E-12&	\textbf{6.3E-13}~/~	1.6E-12&	\textbf{8.6E-13}~/~	2.0E-12&	\textbf{1.1E-12}~/~	2.0E-12\\ \hline
    \multicolumn{1}{|c|}{$\mathbf{|~\norm{S^{sol}}_1-\norm{S^0}_1|/\norm{S^0}_1}$}
    &\textbf{1.1E-5}~/~	1.4E-5&	\textbf{6.2E-6}~/~	9.7E-6&	\textbf{9.7E-6}~/~	1.5E-5&	\textbf{8.6E-5}~/~	9.7E-5\\ \hline
    \multicolumn{1}{|c|}{$\mathbf{\max\{|S^{sol}_{ij}-S^0_{ij}|: S^0_{ij}\neq 0\}}$}
    &\textbf{2.9E-1}~/~	3.5E-1&	\textbf{5.9E-1}~/~	8.0E-1&	\textbf{2.2E-1}~/~	2.4E-1&	\textbf{5.9E-1}~/~	7.4E-1\\ \hline
    \multicolumn{1}{|c|}{$\mathbf{\max\{|S^{sol}_{ij}|: S^0_{ij}=0\}}$}
    &\textbf{0}~/~                0&	\textbf{4.0E-1}~/~	7.2E-1&	\textbf{8.3E-3}~/~	1.1E-2&	\textbf{1.9E-1}~/~	5.5E-1\\ \hline\\
    \cline{2-5}
    &\multicolumn{2}{|c|}{$\mathbf{c_r}$=\textbf{0.1} $\mathbf{c_p}$=\textbf{0.05}}&\multicolumn{2}{|c|}{$\mathbf{c_r}$=\textbf{0.1} $\mathbf{c_p}$=\textbf{0.1}}\\ \cline{2-5}
    &\multicolumn{1}{|c|}{\textbf{NSA}} & \textbf{ASALM}& \textbf{NSA} & \textbf{ASALM} \\ \hline
    \multicolumn{1}{|c|}{\textbf{Error Type}}&\multicolumn{1}{|c|}{\textbf{avg}~/~max}& \textbf{avg}~/~max & \textbf{avg}~/~max & \textbf{avg}~/~max\\ \hline
    \multicolumn{1}{|c|}{$\mathbf{|~\norm{X^{sol}}_*-\norm{X^0}_*|/\norm{X^0}_*}$}
    &\textbf{5.6E-6}~/~	6.4E-6&	\textbf{4.6E-5}~/~	4.9E-5&	\textbf{6.2E-6}~/~	7.1E-6&	\textbf{1.2E-4}~/~	1.4E-4\\ \hline
    \multicolumn{1}{|c|}{$\mathbf{\max\{|\sigma_i-\sigma^0_i|: \sigma^0_i>0\}}$}
    &\textbf{5.7E-2}~/~	6.2E-2&	\textbf{1.2E-1}~/~	1.3E-1&	\textbf{8.8E-2}~/~	1.0E-1&	\textbf{3.0E-1}~/~	3.7E-1\\ \hline
    \multicolumn{1}{|c|}{$\mathbf{\max\{|\sigma_i|: \sigma^0_i=0\}}$}
    &\textbf{6.9E-13}~/~1.5E-12&\textbf{6.2E-13}~/~	9.9E-13&\textbf{6.2E-13}~/~	1.3E-12&\textbf{3.9E-13}~/~	1.0E-12\\ \hline
    \multicolumn{1}{|c|}{$\mathbf{|~\norm{S^{sol}}_1-\norm{S^0}_1|/\norm{S^0}_1}$}
    &\textbf{1.2E-5}~/~	1.6E-5&	\textbf{1.6E-4}~/~	1.7E-4&	\textbf{3.4E-5}~/~	3.7E-5&	\textbf{2.5E-4}~/~	2.7E-4\\ \hline
    \multicolumn{1}{|c|}{$\mathbf{\max\{|S^{sol}_{ij}-S^0_{ij}|: S^0_{ij}\neq 0\}}$}
    &\textbf{1.6E-1}~/~	1.9E-1&	\textbf{6.7E-1}~/~	8.3E-1&	\textbf{1.7E-1}~/~	2.0E-1&	\textbf{7.9E-1}~/~	9.5E-1\\ \hline
    \multicolumn{1}{|c|}{$\mathbf{\max\{|S^{sol}_{ij}|: S^0_{ij}=0\}}$}
    &\textbf{7.0E-3}~/~	1.1E-2&	\textbf{1.5E-1}~/~	2.5E-1&	\textbf{1.3E-2}~/~	1.9E-2&	\textbf{1.2E-1}~/~	2.5E-1\\ \hline
    \end{tabular}
    \label{tab:compare_detail_80dB}
    }
    \end{adjustwidth}
\end{table}

\begin{table}[!htb]
    \begin{adjustwidth}{-2em}{-2em}
    \centering
    \caption{NSA vs ASALM: Additional statistics on solution accuracy for decomposing $D\in\reals^{n\times n}$, $n=1500$, ${\rm SNR}(D)=45dB$}
    \renewcommand{\arraystretch}{1.4}
    {\footnotesize
    \begin{tabular}{cc|c|c|c|}
    \cline{2-5}
    &\multicolumn{2}{|c|}{$\mathbf{c_r}$=\textbf{0.05} $\mathbf{c_p}$=\textbf{0.05}}&\multicolumn{2}{|c|}{$\mathbf{c_r}$=\textbf{0.05} $\mathbf{c_p}$=\textbf{0.1}}\\ \cline{2-5}
    &\multicolumn{1}{|c|}{\textbf{NSA}} & \textbf{ASALM} & \textbf{NSA} & \textbf{ASALM} \\ \hline
    \multicolumn{1}{|c|}{\textbf{Error Type}}&\multicolumn{1}{|c|}{\textbf{avg}~/~max}& \textbf{avg}~/~max & \textbf{avg}~/~max & \textbf{avg}~/~max\\ \hline
    \multicolumn{1}{|c|}{$\mathbf{|~\norm{X^{sol}}_*-\norm{X^0}_*|/\norm{X^0}_*}$}
    &\textbf{1.8E-4}~/~	3.6E-4&	\textbf{1.4E-3}~/~	1.5E-3&	\textbf{2.2E-4}~/~	2.4E-4&	\textbf{2.4E-3}~/~	2.6E-3\\ \hline
    \multicolumn{1}{|c|}{$\mathbf{\max\{|\sigma_i-\sigma^0_i|: \sigma^0_i>0\}}$}
    &\textbf{5.9E-1}~/~	1.8E+0&	\textbf{1.1E+0}~/~	1.5E+0&	\textbf{9.8E-1}~/~	1.1E+0&	\textbf{2.3E+0}~/~	2.6E+0\\ \hline
    \multicolumn{1}{|c|}{$\mathbf{\max\{|\sigma_i|: \sigma^0_i=0\}}$}
    &\textbf{6.4E-13}~/~	1.3E-12&\textbf{3.7E+0}~/~3.8E+0&\textbf{6.1E-13}~/~1.0E-12&	\textbf{4.7E+0}~/~5.5E+0\\ \hline
    \multicolumn{1}{|c|}{$\mathbf{|~\norm{S^{sol}}_1-\norm{S^0}_1|/\norm{S^0}_1}$}
    &\textbf{1.7E-4}~/~	1.9E-4&	\textbf{4.2E-3}~/~	4.3E-3&	\textbf{1.3E-4}~/~	1.3E-4&	\textbf{2.9E-3}~/~	3.6E-3\\ \hline
    \multicolumn{1}{|c|}{$\mathbf{\max\{|S^{sol}_{ij}-S^0_{ij}|: S^0_{ij}\neq 0\}}$}
    &\textbf{1.0E+0}~/~	1.2E+0&	\textbf{3.0E+0}~/~	3.6E+0&	\textbf{1.3E+0}~/~	1.4E+0&	\textbf{3.2E+0}~/~	3.8E+0\\ \hline
    \multicolumn{1}{|c|}{$\mathbf{\max\{|S^{sol}_{ij}|: S^0_{ij}=0\}}$}
    &\textbf{3.6E-1}~/~	4.0E-1&	\textbf{2.2E+0}~/~	2.6E+0&	\textbf{5.3E-1}~/~	6.1E-1&	\textbf{2.3E+0}~/~	3.1E+0\\ \hline\\
    \cline{2-5}
    &\multicolumn{2}{|c|}{$\mathbf{c_r}$=\textbf{0.1} $\mathbf{c_p}$=\textbf{0.05}}&\multicolumn{2}{|c|}{$\mathbf{c_r}$=\textbf{0.1} $\mathbf{c_p}$=\textbf{0.1}}\\ \cline{2-5}
    &\multicolumn{1}{|c|}{\textbf{NSA}} & \textbf{ASALM}& \textbf{NSA} & \textbf{ASALM} \\ \hline
    \multicolumn{1}{|c|}{\textbf{Error Type}}&\multicolumn{1}{|c|}{\textbf{avg}~/~max}& \textbf{avg}~/~max & \textbf{avg}~/~max & \textbf{avg}~/~max\\ \hline
    \multicolumn{1}{|c|}{$\mathbf{|~\norm{X^{sol}}_*-\norm{X^0}_*|/\norm{X^0}_*}$}
    &\textbf{3.7E-4}~/~	6.5E-4&	\textbf{9.7E-5}~/~	1.3E-4&	\textbf{6.7E-4}~/~	6.8E-4&	\textbf{8.4E-4}~/~	9.0E-4\\ \hline
    \multicolumn{1}{|c|}{$\mathbf{\max\{|\sigma_i-\sigma^0_i|: \sigma^0_i>0\}}$}
    &\textbf{1.3E+0}~/~	1.5E+0&	\textbf{1.2E+0}~/~	1.3E+0&	\textbf{2.5E+0}~/~	2.8E+0&	\textbf{1.3E+0}~/~	1.5E+0\\ \hline
    \multicolumn{1}{|c|}{$\mathbf{\max\{|\sigma_i|: \sigma^0_i=0\}}$}
    &\textbf{1.6E-1}~/~	1.6E+0&	\textbf{3.6E+0}~/~	3.7E+0&	\textbf{7.3E-13}~/~	1.7E-12& \textbf{3.2E+0}~/~	3.3E+0\\ \hline
    \multicolumn{1}{|c|}{$\mathbf{|~\norm{S^{sol}}_1-\norm{S^0}_1|/\norm{S^0}_1}$}
    &\textbf{8.1E-4}~/~	3.2E-3&	\textbf{4.7E-3}~/~	4.8E-3&	\textbf{8.9E-4}~/~	9.0E-4&	\textbf{4.4E-3}~/~	4.5E-3\\ \hline
    \multicolumn{1}{|c|}{$\mathbf{\max\{|S^{sol}_{ij}-S^0_{ij}|: S^0_{ij}\neq 0\}}$}
    &\textbf{9.3E-1}~/~	1.1E+0&	\textbf{2.7E+0}~/~	3.3E+0&	\textbf{1.1E+0}~/~	1.2E+0&	\textbf{3.2E+0}~/~	3.5E+0\\ \hline
    \multicolumn{1}{|c|}{$\mathbf{\max\{|S^{sol}_{ij}|: S^0_{ij}=0\}}$}
    &\textbf{5.7E-1}~/~	6.6E-1&	\textbf{1.1E+0}~/~	1.4E+0&	\textbf{7.1E-1}~/~	7.9E-1&	\textbf{1.3E+0}~/~	1.6E+0\\ \hline
    \end{tabular}
    \label{tab:compare_detail_45dB}
    }
    \end{adjustwidth}
\end{table}

\begin{table}[!htb]
    \begin{adjustwidth}{-2em}{-3em}
    \centering
    \caption{NSA: Additional statistics on solution accuracy for decomposing $D\in\reals^{n\times n}$, $n\in\{500, 1000, 1500\}$, ${\rm SNR}(D)=80dB$}
    \renewcommand{\arraystretch}{1.45}
    {\footnotesize
    \begin{tabular}{cc|c|c|c|c|}
    \cline{3-6}
    &&$\mathbf{c_r}$=\textbf{0.05} $\mathbf{c_p}$=\textbf{0.05}&$\mathbf{c_r}$=\textbf{0.05} $\mathbf{c_p}$=\textbf{0.1}&$\mathbf{c_r}$=\textbf{0.1} $\mathbf{c_p}$=\textbf{0.05}&$\mathbf{c_r}$=\textbf{0.1} $\mathbf{c_p}$=\textbf{0.1}\\\hline
    \multicolumn{1}{|c|}{n}& \multicolumn{1}{|c|}{Error Type}& \textbf{avg}~/~max & \textbf{avg}~/~max & \textbf{avg}~/~max & \textbf{avg}~/~max\\ \hline
    \multicolumn{1}{|c|}{\multirow{7}{*}{$\mathbf{500}$}}

    & \multicolumn{1}{|c|}{$\mathbf{\frac{|~\norm{X^{sol}}_*-\norm{X^0}_*|}{\norm{X^0}_*}}$}
    &\textbf{7.2E-6}~/~	1.1E-5&	\textbf{2.0E-5}~/~	2.7E-5&	\textbf{5.6E-6}~/~	8.2E-6&	\textbf{2.1E-5}~/~	3.1E-5\\

    \multicolumn{1}{|c|}{}& \multicolumn{1}{|c|}{$\mathbf{\max\{|\sigma_i-\sigma^0_i|: \sigma^0_i>0\}}$}
    &\textbf{1.7E-2}~/~	2.4E-2&	\textbf{3.4E-2}~/~	5.6E-2&	\textbf{2.1E-2}~/~	2.7E-2&	\textbf{3.2E-2}~/~	3.8E-2\\

    \multicolumn{1}{|c|}{}& \multicolumn{1}{|c|}{$\mathbf{\max\{|\sigma_i|: \sigma^0_i=0\}}$}
    &\textbf{1.6E-13}~/~	2.9E-13&	\textbf{2.0E-13}~/~	5.6E-13&	\textbf{1.1E-13}~/~	2.5E-13&	\textbf{8.6E-14}~/~	1.7E-13\\

    \multicolumn{1}{|c|}{}& \multicolumn{1}{|c|}{$\mathbf{\frac{|~\norm{S^{sol}}_1-\norm{S^0}_1|}{\norm{S^0}_1}}$}
    &\textbf{1.6E-5}~/~	1.7E-5&	\textbf{1.5E-5}~/~	1.8E-5&	\textbf{2.9E-5}~/~	3.2E-5&	\textbf{2.6E-5}~/~	3.0E-5\\

    \multicolumn{1}{|c|}{}& \multicolumn{1}{|c|}{$\mathbf{\max\{|S^{sol}_{ij}-S^0_{ij}|: S^0_{ij}\neq 0\}}$}
    &\textbf{3.2E-1}~/~	4.0E-1&	\textbf{3.0E-1}~/~	4.3E-1&	\textbf{2.6E-1}~/~	3.2E-1&	\textbf{1.8E-1}~/~	2.3E-1\\

    \multicolumn{1}{|c|}{}& \multicolumn{1}{|c|}{$\mathbf{\max\{|S^{sol}_{ij}|: S^0_{ij}=0\}}$}
    &\textbf{9.5E-3}~/~	2.2E-2&	\textbf{1.5E-2}~/~	2.5E-2&	\textbf{1.5E-2}~/~	2.5E-2&	\textbf{1.8E-2}~/~	3.4E-2\\ \hline

    \multicolumn{1}{|c|}{\multirow{7}{*}{$\mathbf{1000}$}}
    & \multicolumn{1}{|c|}{$\mathbf{\frac{|~\norm{X^{sol}}_*-\norm{X^0}_*|}{\norm{X^0}_*}}$}
    &\textbf{5.6E-6}~/~	1.7E-5&	\textbf{6.2E-6}~/~	1.7E-5&	\textbf{6.9E-6}~/~	8.6E-6&	\textbf{1.5E-6}~/~	2.6E-6\\

    \multicolumn{1}{|c|}{}& \multicolumn{1}{|c|}{$\mathbf{\max\{|\sigma_i-\sigma^0_i|: \sigma^0_i>0\}}$}
    &\textbf{1.8E-2}~/~	4.0E-2&	\textbf{3.1E-2}~/~	4.8E-2&	\textbf{5.1E-2}~/~	6.0E-2&	\textbf{5.9E-2}~/~	6.8E-2\\

    \multicolumn{1}{|c|}{}& \multicolumn{1}{|c|}{$\mathbf{\max\{|\sigma_i|: \sigma^0_i=0\}}$}
    &\textbf{3.3E-13}~/~	4.8E-13&	\textbf{3.3E-13}~/~	5.0E-13&	\textbf{2.9E-13}~/~	6.6E-13&	\textbf{2.8E-13}~/~	4.8E-13\\

    \multicolumn{1}{|c|}{}& \multicolumn{1}{|c|}{$\mathbf{\frac{|~\norm{S^{sol}}_1-\norm{S^0}_1|}{\norm{S^0}_1}}$}
    &\textbf{1.1E-5}~/~	1.5E-5&	\textbf{1.7E-5}~/~	1.9E-5&	\textbf{2.8E-5}~/~	3.0E-5&	\textbf{2.9E-5}~/~	3.0E-5\\

    \multicolumn{1}{|c|}{}& \multicolumn{1}{|c|}{$\mathbf{\max\{|S^{sol}_{ij}-S^0_{ij}|: S^0_{ij}\neq 0\}}$}
    &\textbf{2.7E-1}~/~	3.1E-1&	\textbf{3.1E-1}~/~	3.8E-1&	\textbf{2.2E-1}~/~	2.8E-1&	\textbf{1.6E-1}~/~	1.7E-1\\

    \multicolumn{1}{|c|}{}& \multicolumn{1}{|c|}{$\mathbf{\max\{|S^{sol}_{ij}|: S^0_{ij}=0\}}$}
    &\textbf{1.7E-4}~/~	9.7E-4&	\textbf{1.2E-2}~/~	1.7E-2&	\textbf{7.8E-3}~/~	1.2E-2&	\textbf{1.2E-2}~/~	1.5E-2\\ \hline

    \multicolumn{1}{|c|}{\multirow{7}{*}{$\mathbf{1500}$}}
    & \multicolumn{1}{|c|}{$\mathbf{\frac{|~\norm{X^{sol}}_*-\norm{X^0}_*|}{\norm{X^0}_*}}$}
    &\textbf{1.7E-6}~/~	5.2E-6&	\textbf{5.0E-6}~/~	3.2E-5&	\textbf{5.6E-6}~/~	6.4E-6&	\textbf{6.2E-6}~/~	7.1E-6\\

    \multicolumn{1}{|c|}{}& \multicolumn{1}{|c|}{$\mathbf{\max\{|\sigma_i-\sigma^0_i|: \sigma^0_i>0\}}$}
    &\textbf{4.1E-2}~/~	5.2E-2&	\textbf{5.2E-2}~/~	1.8E-1&	\textbf{5.7E-2}~/~	6.2E-2&	\textbf{8.8E-2}~/~	1.0E-1\\

    \multicolumn{1}{|c|}{}& \multicolumn{1}{|c|}{$\mathbf{\max\{|\sigma_i|: \sigma^0_i=0\}}$}
    &\textbf{7.9E-13}~/~	2.2E-12&	\textbf{8.6E-13}~/~	2.0E-12&	\textbf{6.9E-13}~/~	1.5E-12&	\textbf{6.2E-13}~/~	1.3E-12\\

    \multicolumn{1}{|c|}{}& \multicolumn{1}{|c|}{$\mathbf{\frac{|~\norm{S^{sol}}_1-\norm{S^0}_1|}{\norm{S^0}_1}}$}
    &\textbf{1.1E-5}~/~	1.4E-5&	\textbf{9.7E-6}~/~	1.5E-5&	\textbf{1.2E-5}~/~	1.6E-5&	\textbf{3.4E-5}~/~	3.7E-5\\

    \multicolumn{1}{|c|}{}& \multicolumn{1}{|c|}{$\mathbf{\max\{|S^{sol}_{ij}-S^0_{ij}|: S^0_{ij}\neq 0\}}$}
    &\textbf{2.9E-1}~/~	3.5E-1&	\textbf{2.2E-1}~/~	2.4E-1&	\textbf{1.6E-1}~/~	1.9E-1&	\textbf{1.7E-1}~/~	2.0E-1\\

    \multicolumn{1}{|c|}{}& \multicolumn{1}{|c|}{$\mathbf{\max\{|S^{sol}_{ij}|: S^0_{ij}=0\}}$}
    &\textbf{0}~/~	0&	\textbf{8.3E-3}~/~	1.1E-2&	\textbf{7.0E-3}~/~	1.1E-2&	\textbf{1.3E-2}~/~	1.9E-2\\ \hline
    \end{tabular}
    \label{tab:self_detail_80dB}
    }
    \end{adjustwidth}
\end{table}

\begin{table}[!htb]
    \begin{adjustwidth}{-2em}{-3em}
    \centering
    \caption{NSA: Additional statistics on solution accuracy for decomposing $D\in\reals^{n\times n}$, $n\in\{500, 1000, 1500\}$, ${\rm SNR}(D)=45dB$}
    \renewcommand{\arraystretch}{1.45}
    {\footnotesize
    \begin{tabular}{cc|c|c|c|c|}
    \cline{3-6}
    &&$\mathbf{c_r}$=\textbf{0.05} $\mathbf{c_p}$=\textbf{0.05}&$\mathbf{c_r}$=\textbf{0.05} $\mathbf{c_p}$=\textbf{0.1}&$\mathbf{c_r}$=\textbf{0.1} $\mathbf{c_p}$=\textbf{0.05}&$\mathbf{c_r}$=\textbf{0.1} $\mathbf{c_p}$=\textbf{0.1}\\\hline
    \multicolumn{1}{|c|}{n}& \multicolumn{1}{|c|}{Error Type}& \textbf{avg}~/~max & \textbf{avg}~/~max & \textbf{avg}~/~max & \textbf{avg}~/~max\\ \hline
    \multicolumn{1}{|c|}{\multirow{7}{*}{$\mathbf{500}$}}
    & \multicolumn{1}{|c|}{$\mathbf{\frac{|~\norm{X^{sol}}_*-\norm{X^0}_*|}{\norm{X^0}_*}}$}
    &\textbf{6.0E-4}~/~	9.3E-4&	\textbf{5.5E-4}~/~	6.2E-4&	\textbf{7.4E-4}~/~	8.8E-4&	\textbf{1.0E-3}~/~	1.3E-3\\

    \multicolumn{1}{|c|}{}& \multicolumn{1}{|c|}{$\mathbf{\max\{|\sigma_i-\sigma^0_i|: \sigma^0_i>0\}}$}
    &\textbf{5.1E-1}~/~	7.8E-1&	\textbf{5.4E-1}~/~	7.7E-1&	\textbf{8.2E-1}~/~	8.9E-1&	\textbf{9.2E-1}~/~	1.2E+0\\

    \multicolumn{1}{|c|}{}& \multicolumn{1}{|c|}{$\mathbf{\max\{|\sigma_i|: \sigma^0_i=0\}}$}
    &\textbf{1.7E-13}~/~2.7E-13&\textbf{1.6E-13}~/~	3.0E-13&\textbf{1.0E-13}~/~	2.1E-13&\textbf{1.1E-1}~/~	6.0E-1\\

    \multicolumn{1}{|c|}{}& \multicolumn{1}{|c|}{$\mathbf{\frac{|~\norm{S^{sol}}_1-\norm{S^0}_1|}{\norm{S^0}_1}}$}
    &\textbf{3.0E-4}~/~	3.4E-4&	\textbf{2.1E-4}~/~	2.9E-4&	\textbf{3.0E-4}~/~	1.2E-3&	\textbf{6.4E-4}~/~	1.1E-3\\

    \multicolumn{1}{|c|}{}& \multicolumn{1}{|c|}{$\mathbf{\max\{|S^{sol}_{ij}-S^0_{ij}|: S^0_{ij}\neq 0\}}$}
    &\textbf{1.6E+0}~/~	1.9E+0&	\textbf{1.4E+0}~/~	1.8E+0&	\textbf{1.2E+0}~/~	1.6E+0&	\textbf{1.0E+0}~/~	1.3E+0\\

    \multicolumn{1}{|c|}{}& \multicolumn{1}{|c|}{$\mathbf{\max\{|S^{sol}_{ij}|: S^0_{ij}=0\}}$}
    &\textbf{2.3E-1}~/~	2.9E-1&	\textbf{4.0E-1}~/~	4.7E-1&	\textbf{3.5E-1}~/~	4.6E-1&	\textbf{5.4E-1}~/~	6.1E-1\\ \hline

    \multicolumn{1}{|c|}{\multirow{7}{*}{$\mathbf{1000}$}}
    & \multicolumn{1}{|c|}{$\mathbf{\frac{|~\norm{X^{sol}}_*-\norm{X^0}_*|}{\norm{X^0}_*}}$}
    &\textbf{2.8E-4}~/~	3.1E-4&	\textbf{4.4E-4}~/~	7.5E-4&	\textbf{5.6E-4}~/~	8.0E-4&	\textbf{7.4E-4}~/~	8.4E-4\\

    \multicolumn{1}{|c|}{}& \multicolumn{1}{|c|}{$\mathbf{\max\{|\sigma_i-\sigma^0_i|: \sigma^0_i>0\}}$}
    &\textbf{5.2E-1}~/~	6.2E-1&	\textbf{8.6E-1}~/~	1.2E+0&	\textbf{1.7E+0}~/~	1.9E+0&	\textbf{1.8E+0}~/~	1.9E+0\\

    \multicolumn{1}{|c|}{}& \multicolumn{1}{|c|}{$\mathbf{\max\{|\sigma_i|: \sigma^0_i=0\}}$}
    &\textbf{2.5E-13}~/~5.3E-13&\textbf{4.3E-13}~/~	9.0E-13&\textbf{2.0E-1}~/~	2.0E+0&	\textbf{6.3E-1}~/~	3.9E+0\\

    \multicolumn{1}{|c|}{}& \multicolumn{1}{|c|}{$\mathbf{\frac{|~\norm{S^{sol}}_1-\norm{S^0}_1|}{\norm{S^0}_1}}$}
    &\textbf{2.2E-4}~/~	2.3E-4&	\textbf{1.4E-4}~/~	1.7E-4&	\textbf{5.5E-4}~/~	3.7E-3&	\textbf{1.1E-3}~/~	2.5E-3\\

    \multicolumn{1}{|c|}{}& \multicolumn{1}{|c|}{$\mathbf{\max\{|S^{sol}_{ij}-S^0_{ij}|: S^0_{ij}\neq 0\}}$}
    &\textbf{1.3E+0}~/~	1.5E+0&	\textbf{1.5E+0}~/~	1.8E+0&	\textbf{1.1E+0}~/~	1.3E+0&	\textbf{9.6E-1}~/~	1.1E+0\\

    \multicolumn{1}{|c|}{}& \multicolumn{1}{|c|}{$\mathbf{\max\{|S^{sol}_{ij}|: S^0_{ij}=0\}}$}
    &\textbf{2.7E-1}~/~	3.2E-1&	\textbf{4.6E-1}~/~	5.2E-1&	\textbf{4.6E-1}~/~	5.1E-1&	\textbf{6.4E-1}~/~	6.7E-1\\ \hline

    \multicolumn{1}{|c|}{\multirow{7}{*}{$\mathbf{1500}$}}
    & \multicolumn{1}{|c|}{$\mathbf{\frac{|~\norm{X^{sol}}_*-\norm{X^0}_*|}{\norm{X^0}_*}}$}
    &\textbf{1.8E-4}~/~	3.6E-4&	\textbf{2.2E-4}~/~	2.4E-4&	\textbf{3.7E-4}~/~	6.5E-4&	\textbf{6.7E-4}~/~	6.8E-4\\

    \multicolumn{1}{|c|}{}& \multicolumn{1}{|c|}{$\mathbf{\max\{|\sigma_i-\sigma^0_i|: \sigma^0_i>0\}}$}
    &\textbf{5.9E-1}~/~	1.8E+0&	\textbf{9.8E-1}~/~	1.1E+0&	\textbf{1.3E+0}~/~	1.5E+0&	\textbf{2.5E+0}~/~	2.8E+0\\

    \multicolumn{1}{|c|}{}& \multicolumn{1}{|c|}{$\mathbf{\max\{|\sigma_i|: \sigma^0_i=0\}}$}
    &\textbf{6.4E-13}~/~1.3E-12&\textbf{6.1E-13}~/~	1.0E-12&\textbf{1.6E-1}~/~	1.6E+0&	\textbf{7.3E-13}~/~	1.7E-12\\

    \multicolumn{1}{|c|}{}& \multicolumn{1}{|c|}{$\mathbf{\frac{|~\norm{S^{sol}}_1-\norm{S^0}_1|}{\norm{S^0}_1}}$}
    &\textbf{1.7E-4}~/~	1.9E-4&	\textbf{1.3E-4}~/~	1.3E-4&	\textbf{8.1E-4}~/~	3.2E-3&	\textbf{8.9E-4}~/~	9.0E-4\\

    \multicolumn{1}{|c|}{}& \multicolumn{1}{|c|}{$\mathbf{\max\{|S^{sol}_{ij}-S^0_{ij}|: S^0_{ij}\neq 0\}}$}
    &\textbf{1.0E+0}~/~	1.2E+0&	\textbf{1.3E+0}~/~	1.4E+0&	\textbf{9.3E-1}~/~	1.1E+0&	\textbf{1.1E+0}~/~	1.2E+0\\

    \multicolumn{1}{|c|}{}& \multicolumn{1}{|c|}{$\mathbf{\max\{|S^{sol}_{ij}|: S^0_{ij}=0\}}$}
    &\textbf{3.6E-1}~/~	4.0E-1&	\textbf{5.3E-1}~/~	6.1E-1&	\textbf{5.7E-1}~/~	6.6E-1&	\textbf{7.1E-1}~/~	7.9E-1\\ \hline
    \end{tabular}
    \label{tab:self_detail_45dB}
    }
    \end{adjustwidth}
\end{table}
\end{document}